\documentclass[11pt,a4paper]{article}
\usepackage[left=2.5cm,right=2.5cm,top=3cm,bottom=3cm]{geometry}
\usepackage[T1]{fontenc}
\usepackage{lmodern}
\usepackage{amsmath,amssymb,amsthm,bm,mathrsfs}
\usepackage{graphicx}
\usepackage{booktabs,multirow}
\usepackage{caption}
\usepackage{placeins}
\usepackage{authblk}
\usepackage{titlesec}
\usepackage{setspace}
\usepackage{microtype}
\usepackage{xcolor}
\usepackage[hidelinks]{hyperref}

\titleformat{\section}
  {\fontsize{13pt}{15pt}\selectfont\bfseries}{\thesection}{1em}{}
\titleformat{\subsection}
  {\fontsize{12pt}{15pt}\selectfont\bfseries}{\thesubsection}{1em}{}
\numberwithin{equation}{section}
\newtheorem{theorem}{Theorem}[section]
\newtheorem{proposition}[theorem]{Proposition}
\newtheorem{corollary}[theorem]{Corollary}

\title{Asymptotic-Preserving Tensor Neural Networks for 3D High-Field Kinetic Equations with Applications to Uncertainty Quantification}
\author[1,2,3]{Shi Jin}
\author[1]{Tian-ai Zhang\thanks{Corresponding author: \href{mailto:zhangtianai@sjtu.edu.cn}{zhangtianai@sjtu.edu.cn}.}}
\affil[1]{School of Mathematical Sciences, Shanghai Jiao Tong University, Shanghai, China}
\affil[2]{Institute of Natural Sciences, MOE-LSC, Shanghai Jiao Tong University, Shanghai, China}
\affil[3]{Shanghai Artificial Intelligence Laboratory, Shanghai, China}

\date{}
\hypersetup{
  pdftitle={APTNN for 3D High-Field Kinetic Equations},
  pdfauthor={Shi Jin and Tian-ai Zhang}
}

\begin{document}
\maketitle

\begin{abstract}
Kinetic equations in the high-field regime describe particle transport dominated by electric-field acceleration and collisional relaxation. These equations are widely used in simulations of plasmas and semiconductor devices. Typical models include the Vlasov--Poisson--Fokker--Planck (VPFP) and Boltzmann--Poisson systems. Their high dimensionality and multiscale nature pose two of the major challenges to accurate and efficient numerical solution. To address these challenges, we propose an asymptotic-preserving tensor neural network (APTNN) method. Specifically, to reduce high-dimensional computational costs and memory requirements, we adopt the canonical polyadic (CP) decomposition structure of TNNs to approximate their solutions. We show that the equation residuals retain this decomposition form, allowing high-dimensional loss integrals to be evaluated through combinations of one-dimensional integrals. The resulting computational and storage costs are reduced from exponential to polynomial growth with dimension. To improve the prediction accuracy in the high-field limit, we enhance the baseline AP loss by incorporating the first moment equation. The added loss term directly constrains the limiting flux closure to reduce its approximation error in the high-field regime. The enhanced AP loss can still be evaluated dimension-wise. Numerical results show that APTNN improves prediction accuracy in the high-field regime with significantly reduced computational costs and memory usage. For a 3D3V VPFP problem with multiple random inputs, training takes less than $30$ minutes on a single GPU.
\end{abstract}


\section{Introduction}

Kinetic equations provide a statistical description of motion and interactions in a collection of particles \cite{dimarco2014numerical}. They are fundamental to modeling plasma dynamics and semiconductor transport \cite{dawson1983particle, jacoboni1983monte}. Typical models include the Vlasov--Poisson--Fokker--Planck (VPFP) system for plasmas \cite{JinWang2011,GuZhang2026} and the Boltzmann--Poisson system for semiconductors \cite{JinWang2013}, both of which incorporate a self-consistent electric field. Of great physical interest is the high-field regime, in which electric-field acceleration and collisional relaxation jointly dominate particle transport~\cite{NietoPoupaudSoler,CercignaniGambaLevermore,CoughlinHu2022}. This regime is important for modeling the Brownian motion of charged particles in a thermal bath and hot-electron transport in miniaturized semiconductor devices~\cite{JinWang2011,JinWang2013}.

Numerically solving kinetic equations in the high-field regime presents two major challenges: high dimensionality and multiple temporal and spatial scales. First, the distribution function is defined over phase space and time, involving seven variables in three-dimensional problems. For full tensor-product discretizations, the number of grid points grows exponentially with dimension, leading to the curse of dimensionality~\cite{LiEtAl2023}. Moreover, uncertainty in initial data, boundary conditions or model parameters introduces additional random variables, which further increase the computational cost of uncertainty quantification (UQ)~\cite{JinShaoWang}. Second, the multiscale nature of these equations presents another challenge. As the characteristic scales become small, numerically resolving them can be prohibitively expensive. Therefore, an effective numerical scheme should reduce high-dimensional computational costs and capture the correct macroscopic limit without resolving the small scales. 

In recent years, neural network (NN) methods have emerged as promising approaches to solving high-dimensional partial differential equations (PDEs)~\cite{HanJentzenE2018,EYu2018,WangXinZhang2022}. These methods leverage the powerful approximation capabilities of NNs to predict PDE solutions. For example, physics-informed neural networks (PINNs)~\cite{RaissiEtAl2019} compute derivatives using automatic differentiation (AD) and construct a loss from mean squared PDE residuals. Residuals of the initial and boundary conditions (IBCs) are also included. To improve computational efficiency in high dimensions, tensor neural networks (TNNs) employ a low-rank decomposition idea to modify NN architectures \cite{WangJinXie2024}. They feed each coordinate into a separate subnetwork and combine the outputs as sums of products to represent multivariate functions. This separable structure is also adopted in separable physics-informed neural networks (SPINNs) to reduce repeated network computations in function and derivative evaluation~\cite{ChoEtAl2023}. Recent studies have applied this approach to kinetic equations to reduce computational costs~\cite{OhEtAl2025,ChenDimarcoPareschi2025}. However, the high dimensionality of kinetic equations makes loss evaluation on full tensor-product grids computationally expensive and memory intensive. In particular, GPU memory requirements grow exponentially with dimension and become a primary constraint on problem size. This bottleneck becomes more severe in UQ problems with additional random inputs.

Besides addressing the challenges of high dimensionality, NNs also need to accurately capture the multiscale behavior of kinetic equations. The asymptotic-preserving (AP) idea has been incorporated into neural solvers to preserve the correct macroscopic limit, leading to asymptotic-preserving neural networks (APNNs)~\cite{JinMaWu2023,LuWangXu2022,BaiChenWu2025}. For the high-field regime, APNN can be constructed by combining the macroscopic mass conservation law with the kinetic equations in the loss \cite{JinMaZhang2024,ZhangJin2025}. However, the correct limiting behavior alone does not guarantee prediction accuracy in the limit. The baseline AP loss lacks a direct constraint on the limiting flux closure for the continuity equation. A misfit in this closure may reduce the predicted accuracy for the macroscopic evolution. 

In this paper, we aim to develop an efficient and accurate NN method for high-dimensional kinetic equations in the high-field regime. First, we adopt TNNs approximation in variable-separation form. For their training loss, we show that a dimension-wise integration strategy can be introduced to express the high-dimensional losses as combinations of one-dimensional integrals. This avoids the full-grid loss evaluation used in previous TNN-based kinetic solvers~\cite{OhEtAl2025,ChenDimarcoPareschi2025}, reducing computational costs and memory usage from exponential to polynomial growth with dimension. Second, we enhance the AP loss by incorporating the first velocity-moment equation in the loss term. The added term directly constrains the flux closure  and improves predictions of the macroscopic evolution in the high-field limit. The enhanced loss remains separable and allows efficient dimension-wise evaluation. As a result, we developed an asymptotic-preserving tensor neural network (APTNN) method. Numerical results show that APTNN improves prediction accuracy in the high-field limit while significantly reducing computational costs and memory usage. A multidimensional UQ example for a 3D3V problem requires less than $30$ minutes of training on a single GPU.

The remainder of this paper is organized as follows. Sect.~2 introduces the kinetic equations in the high-field regime and their high-field limits. Sect.~3 presents the APTNN method, starting with the separable TNN structure and its approximation properties. It then constructs an AP loss with an enhanced constraint on the limiting flux closure and derives a dimension-wise formula for loss evaluation. Numerical results in Sect.~4 demonstrate the accuracy and efficiency of the proposed method across different dimensions and scale parameters. Multi-dimensional UQ examples are also included. Sect.~5 concludes the paper and discusses future research directions. 

\section{Model equations}
This section introduces kinetic models in the high-field regime, where both electric force and collisions dominate the transport of charged particles. Typical applications include the simulation of plasmas and semiconductor devices. Let $f=f(t,\mathbf{x},\mathbf{v})\geq0$ denote the distribution function of electrons  at time $t\in[0,T]$ in the phase space $\Omega_\mathbf{x}\times\mathbb{R}^{d}$, where $\mathbf{x}\in\Omega_\mathbf{x}\subset\mathbb{R}^{d}$, $\mathbf{v}\in\mathbb{R}^{d}$ and $d\in\{1,2,3\}$. Its time evolution under a self-consistent electrostatic potential $\phi$ is governed by
\begin{equation}\label{eq:general_kinetic_model}
    \begin{cases}
        \displaystyle
        \partial_t f
        +\mathbf{v}\cdot\nabla_{\mathbf{x}}f
        +\dfrac{1}{\varepsilon}
        \mathbf{E}\cdot\nabla_{\mathbf{v}}f
        =
        \dfrac{1}{\varepsilon}\mathcal{Q}(f),
        \\[2mm]
        \displaystyle
        -\Delta_{\mathbf{x}}\phi
        =\rho-h,
        \qquad
        \mathbf{E}=-\nabla_{\mathbf{x}}\phi.
    \end{cases}
\end{equation}
Here, $\mathbf{E}$ is the electric field and $\mathcal{Q}$ is a collision operator satisfying mass conservation 
\[\int_{\mathbb{R}^{d}}\mathcal{Q}(f)\,\mathrm{d}\mathbf{v}=0.
\]
In the second Poisson equation, 
\[
\rho(t,\mathbf{x})=\int_{\mathbb{R}^{d}}f(t,\mathbf{x},\mathbf{v})\,\mathrm{d}\mathbf{v}
\]
is the electron density and $h=h(\mathbf{x})$ is the given background charge density. The system satisfies the global charge-neutrality condition
\begin{equation}\label{eq:global_charge_neutrality}
    \int_{\Omega_\mathbf{x}}\int_{\mathbb{R}^{d}}
    f^{0}(\mathbf{x},\mathbf{v})
    \,\mathrm{d}\mathbf{v}\,\mathrm{d}\mathbf{x}=\int_{\Omega_\mathbf{x}}h(\mathbf{x})\,\mathrm{d}\mathbf{x},
\end{equation}
where $f^{0}=f(0,\mathbf{x},\mathbf{v})$ is the initial distribution and $h$ is set to $1$ in this work. The scale parameter $\varepsilon>0$ denotes the ratio of the electron mean free path to the characteristic length of the system. When $\varepsilon=O(1)$, the system is in the kinetic regime, whereas $\varepsilon\ll1$ corresponds to the high-field regime. The asymptotic limit $\varepsilon\to0$ is referred to as the high-field limit.

We consider two collision operators, namely the Fokker--Planck operator for the VPFP system and the Boltzmann operator for the semiconductor Boltzmann equation. For the VPFP system in electrostatic plasmas, the Fokker--Planck operator models Coulomb interactions between electrons and a surrounding thermal bath \cite{JinWang2011} and is given by 
\begin{equation}\label{eq:fokker_planck_operator}
    \mathcal{Q}_{\mathrm{FP}}(f)
    =
    \nabla_{\mathbf{v}}\cdot
    \left(
        \mathbf{v}f+\nabla_{\mathbf{v}}f
    \right).
\end{equation}
Taking $\mathcal{Q}=\mathcal{Q}_{\mathrm{FP}}$ in \eqref{eq:general_kinetic_model} and combining it with the electric force term, the Vlasov equation can be rewritten as
\begin{equation}\label{eq:vpfp_system}
    \partial_t f
    +\mathbf{v}\cdot\nabla_{\mathbf{x}}f
    =
    \frac{1}{\varepsilon}\mathcal{L}_{\mathbf{E}}f,
    \quad \text{where }
    \mathcal{L}_{\mathbf{E}}f
    := \nabla_{\mathbf{v}}\cdot
    \left[
        (\mathbf{v}-\mathbf{E})f
        +\nabla_{\mathbf{v}}f
    \right].
\end{equation}
The combined operator $\mathcal{L}_{\mathbf{E}}$ remains of Fokker--Planck type. Its null space is 
\[
\mathcal{N}(\mathcal{L}_{\mathbf{E}}) =\operatorname{Span}\{M_{\mathbf{E}}\} =\{f=\rho M_{\mathbf{E}}\},
\]
where 
\[
M_{\mathbf{E}}(\mathbf{v})=(2\pi)^{-d/2}\exp(-|\mathbf{v}-\mathbf{E}|^{2}/2)
\]
is the shifted Maxwellian that depends on $(t,\mathbf{x})$ through $\mathbf{E}$. 

For the semiconductor Boltzmann equation, one considers the low density regime, in which the collision operator takes the form 
\begin{equation}\label{eq:semiconductor_collision_operator}
    \mathcal{Q}_{\mathrm{semiB}}(f)(\mathbf{v})
    =
    \int_{\mathbb{R}^{d}}
    \Psi(\mathbf{v},\mathbf{v}')
    \left[M(\mathbf{v})f(\mathbf{v}')-M(\mathbf{v}')f(\mathbf{v})
    \right]
    \,\mathrm{d}\mathbf{v}',
\end{equation}
which is called the non-degenerate case \cite{JinWang2013}. Here, $\Psi(\mathbf{v},\mathbf{v}')=\Psi(\mathbf{v}',\mathbf{v})$ is the symmetric scattering kernel and 
\[
M(\mathbf{v}) =(2\pi)^{-d/2}\exp(-|\mathbf{v}|^{2}/2)
\]
is the standard Maxwellian. This operator characterizes the interactions of electrons with the semiconductor crystal lattice and with one another. In contrast to the VPFP system, the high-field equilibrium of the semiconductor Boltzmann equation is not always available in explicit form. In this work, we consider the simplest time-relaxation approximation $\Psi\equiv1$, for which
\begin{equation}\label{eq:bgk_collision_operator}
    \mathcal{Q}_{\mathrm{BGK}}(f)
    =
    \rho M-f.
\end{equation}
With this collision, \eqref{eq:general_kinetic_model} defines the BGK--Poisson system. 

We next show the common high-field limit of the VPFP and BGK--Poisson systems. As $\varepsilon\to0$, $f$ approaches a local high-field equilibrium. For the VPFP system, $\mathcal{N}(\mathcal{L}_{\mathbf{E}})=\operatorname{Span}\{M_{\mathbf{E}}\}$ implies $f\to\rho M_{\mathbf{E}}$, whereas for the BGK--Poisson system, $f\to\rho F_{\mathbf{E}}$, where $F_{\mathbf{E}}$ satisfies
\begin{equation}\label{eq:bgk_high_field_equilibrium}
    \int_{\mathbb{R}^{d}}
    F_{\mathbf{E}}\,\mathrm{d}\mathbf{v}=1,
    \quad
    \mathbf{E}\cdot\nabla_{\mathbf{v}}F_{\mathbf{E}}
    =M-F_{\mathbf{E}},
    \quad
    F_{\mathbf{E}}\geq0.
\end{equation}
Since both $\mathcal{Q}_{\mathrm{FP}}$ and $\mathcal{Q}_{\mathrm{BGK}}$ conserve mass, integrating \eqref{eq:general_kinetic_model} over velocity space yields the continuity equation
\begin{equation}\label{eq:macroscopic_mass_conservation}
    \partial_t\rho
    +\nabla_{\mathbf{x}}\cdot\mathbf{j}=0,
\end{equation}
where $\mathbf{j}=\int_{\mathbb{R}^{d}}\mathbf{v}f\,\mathrm{d}\mathbf{v}$ is the flux. Because both collision operators also satisfy the same first velocity-moment identity
\begin{equation}\label{eq:collision_first_moment}
    \int_{\mathbb{R}^{d}}\mathbf{v}\mathcal{Q}(f)\,\mathrm{d}\mathbf{v}
    =-\mathbf{j},
\end{equation}
multiplying \eqref{eq:general_kinetic_model} by $\mathbf{v}$ and integrating over $\mathbb{R}^{d}$ gives
\begin{equation}\label{eq:macroscopic_momentum_equation}
    \varepsilon
    \left(\partial_t\mathbf{j}
        +\nabla_{\mathbf{x}}\cdot\boldsymbol{\Pi}
    \right)
    +\mathbf{j}-\rho\mathbf{E}
    =0,
\end{equation}
where $\boldsymbol{\Pi}=\int_{\mathbb{R}^{d}}\mathbf{v}\otimes\mathbf{v}f\,\mathrm{d}\mathbf{v}$ is the momentum flux tensor. Letting $\varepsilon\to0$ in this equation yields $\mathbf{j}=\rho\mathbf{E}$. Substituting this relation into \eqref{eq:macroscopic_mass_conservation}, together with the Poisson equation, one obtains the high-field limit system
\begin{equation}\label{eq:high_field_limit}
    \begin{cases}
        \displaystyle
        \partial_t\rho
        +
        \nabla_{\mathbf{x}}\cdot\left(
            \rho\mathbf{E}
        \right)=0, \\[2mm]
        \displaystyle
        -\Delta_{\mathbf{x}}\phi=\rho-h,
        \qquad
        \mathbf{E}=-\nabla_{\mathbf{x}}\phi.
    \end{cases}
\end{equation}
This formal analysis can be made rigorous. For the VPFP system, the high-field limit was first established in one dimension \cite{NietoPoupaudSoler} and later extended to multiple dimensions \cite{GoudonEtAl}. This result for the BGK--Poisson system was established in \cite{CercignaniGambaLevermore}. For more general Boltzmann operators, this topic remains an active area of research.

\section{Methodology} 

This section presents the construction of the APTNN method. We first introduce the TNN architecture and establish an approximation result on the unbounded velocity domain. Then an enhanced AP loss function is formulated to reduce the approximation error in the limiting flux closure. Finally, we derive a dimension-wise formula for loss evaluation and analyze its computational cost and storage requirements.

\subsection{The tensor neural network architecture} 
Standard PINNs can incur substantially higher computational costs as the input dimension $D$ grows. This mainly stems from using coordinate-based network architectures that map input coordinates to the corresponding solution values. For \eqref{eq:general_kinetic_model}, the main high-dimensional solution is the distribution $f(t,\mathbf{x},\mathbf{v})$ with $D_f=1+d_x+d_v$ variables, where $d_x$ and $d_v$ denote the spatial and velocity dimensions respectively. Introducing $d_z$ random variables in UQ problems further increases the input dimension to $D_f=1+d_x+d_v+d_z$. In the PINN framework, $f$ is approximated by a network $f_\theta$, where $\theta$ collects the trainable parameters. If $N$ collocation points are used in each coordinate, $f_\theta$ and its derivatives are evaluated at $N^{D_f}$ tensor-product points. 

A well-established approach to tackling high-dimensional challenges is canonical polyadic (CP) tensor decomposition~\cite{kolda2009tensor,bachmayr2023low}. Its key idea is to represent a high-dimensional tensor as a finite sum of rank-one tensor products. Such separated representations are also suitable for approximating solutions of kinetic equations, since these solutions often evolve near a low-dimensional manifold and exhibit a low-rank structure~\cite{DektorRodgersVenturi2021,ZhangHu2026}. Thus we can approximate $f(t,\mathbf{x},\mathbf{v})$ in CP decomposition form
\begin{equation}\label{eq:tnn_distribution}
    f(t,\mathbf{x},\mathbf{v})
    \approx
    \sum_{r=1}^{r_f} c_r^f
    \xi_r^t(t)
    \prod_{i=1}^{d_x}\xi_r^{x_i}(x_i)
    \prod_{j=1}^{d_v}\xi_r^{v_j}(v_j),
\end{equation}
where $r_f$ is the separation rank and $\{c_r^f\}_{r=1}^{r_f}$ are the weights. For UQ problems, each rank-one term additionally contains $\prod_{m=1}^{d_z}\xi_r^{z_m}(z_m)$. 

\begin{figure}[tbp]
    \centering
    \includegraphics[width=\textwidth]{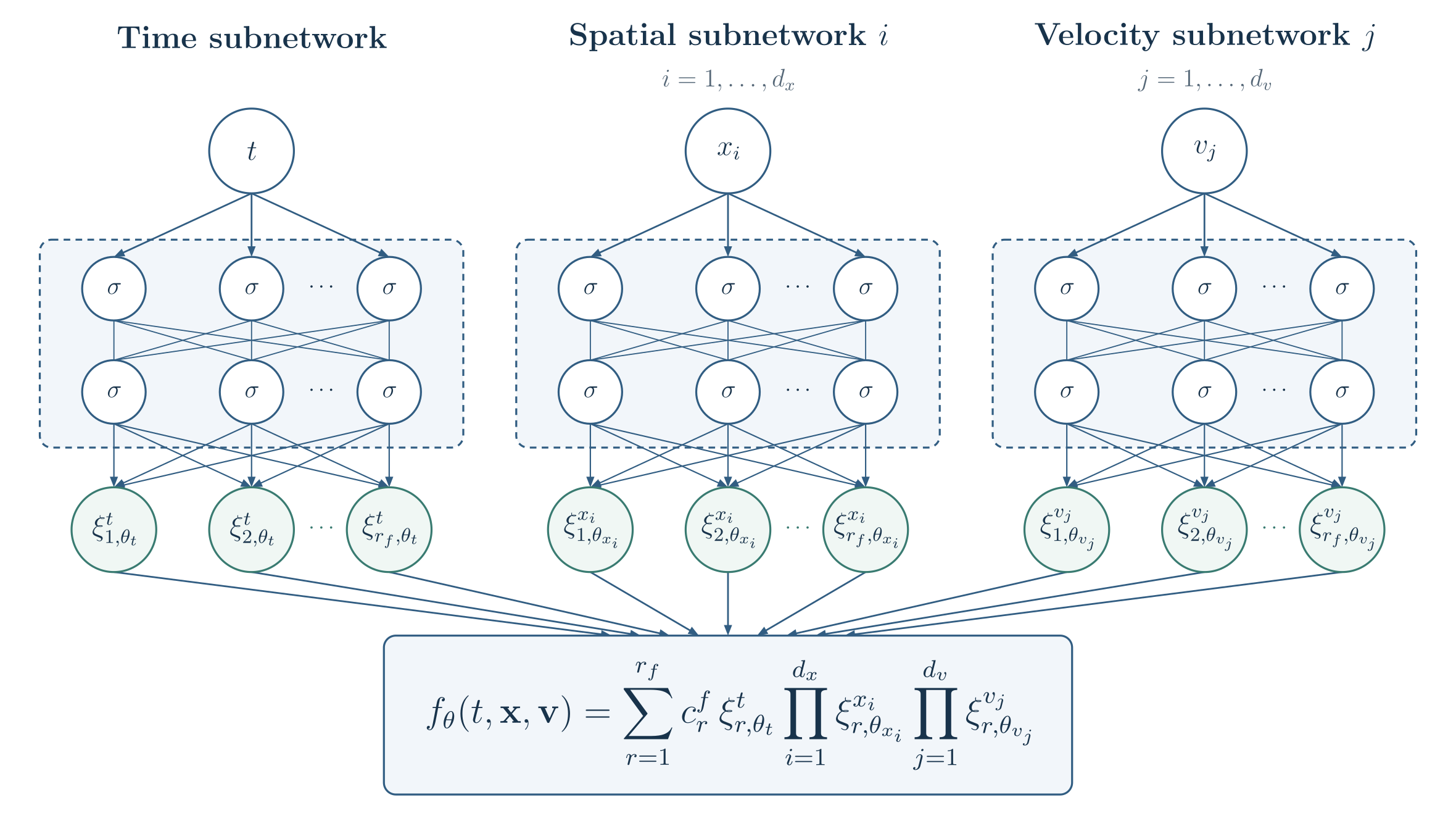}
    \caption{Tensor neural network architecture for $f_\theta$. Each scalar coordinate is assigned an independent subnetwork. Only one spatial and one velocity subnetwork are shown. Here, $\sigma$ denotes the activation function, and $\xi^s_{r,\theta_s}$ denotes the $r$th output component $\xi_r^s(s;\theta_s)$ of the subnetwork for coordinate $s$. To obtain the final output, subnetwork outputs with the same rank index are multiplied and combined using the trainable weights $c_r^f$.}
    \label{fig:tnn_architecture}
\end{figure}

We use NNs to parameterize the CP representation of $f$ in \eqref{eq:tnn_distribution}, which is referred to as a tensor neural network (TNN) \cite{WangJinXie2024, ChoEtAl2023}. The TNN comprises $D_f$ independent feedforward subnetworks, each of which takes a scalar coordinate as input. Let $s$ denote any coordinate, namely, $s=t$, $s=x_i$ for $i=1,\ldots,d_x$, or $s=v_j$ for $j=1,\ldots,d_v$. The corresponding subnetwork maps $s$ to an $r_f$-dimensional output vector
\begin{equation}\label{eq:subnet_vec}
    \boldsymbol{\xi}^{s}(s;\theta_s)
    =\big(
        \xi^{s}_1(s;\theta_s),\ldots,\xi^{s}_{r_f}(s;\theta_s)
    \big)^{\top}.
\end{equation}
Each subnetwork is a feedforward neural network (FNN) with $L$ hidden layers
\begin{equation}\label{eq:subnet_ff}
    \boldsymbol{\xi}^{s}(s;\theta_s)
    =
    \left(
        \boldsymbol{F}_{s}^{(L+1)}
        \circ\boldsymbol{F}_{s}^{(L)}
        \circ\cdots\circ
        \boldsymbol{F}_{s}^{(1)}
    \right)(s).
\end{equation}
For $l=1,\ldots,L$, the hidden layers are defined by $\boldsymbol{F}_{s}^{(l)}(\mathbf{h})=\sigma(\mathbf{W}_{s}^{(l)}\mathbf{h}+\mathbf{b}_{s}^{(l)})$, where $\sigma$ is the activation function. The output layer $\boldsymbol{F}_{s}^{(L+1)}$ is affine. The weight matrices and bias satisfy $\mathbf{W}_{s}^{(l)}\in\mathbb{R}^{N_l\times N_{l-1}}$ and $\mathbf{b}_{s}^{(l)}\in\mathbb{R}^{N_l}$, where $[N_0,N_1,\ldots,N_L,N_{L+1}]$ denotes the layer widths with $N_0=1$ and $N_{L+1}=r_f$. The parameters of each subnetwork are collected in $\theta_s=\{\mathbf{W}_{s}^{(l)}, \mathbf{b}_{s}^{(l)}\}_{l=1}^{L+1}$. The coefficients $\{c_r^f\}_{r=1}^{r_f}$ in \eqref{eq:tnn_distribution} are also trainable parameters. The final output is obtained by multiplying the subnetwork outputs componentwise and taking a weighted sum
\begin{equation}\label{eq:tnn_f_theta}
    f_\theta(t,\mathbf{x},\mathbf{v})
    =
    \sum_{r=1}^{r_f} c_r^f
    \xi_r^t(t;\theta_t)
    \prod_{i=1}^{d_x}\xi_r^{x_i}(x_i;\theta_{x_i})
    \prod_{j=1}^{d_v}\xi_r^{v_j}(v_j;\theta_{v_j}),
\end{equation}
where $\theta$ collects all subnetwork parameters and the trainable coefficients $\{c_r^f\}_{r=1}^{r_f}$. Figure~\ref{fig:tnn_architecture} shows the TNN architecture for $f_\theta$. The potential $\phi$ is also approximated by an independent TNN
\begin{equation}\label{eq:tnn_potential}
    \phi_\theta(t,\mathbf{x})
    =
    \sum_{r=1}^{r_\phi} c_r^\phi
    \eta_r^t(t)
    \prod_{i=1}^{d_x}\eta_r^{x_i}(x_i),
\end{equation}
where $r_\phi$ is the separation rank and $\{c_r^\phi\}_{r=1}^{r_\phi}$ are trainable coefficients.

\begin{samepage}
We establish a TNN approximation theorem in Sobolev spaces for the distribution $f$ on the whole velocity space $\mathbb{R}^3$. We first recall the approximation result in $H^m(\Omega_1\times\cdots\times\Omega_d)$ on bounded domains \cite[Theorem~1]{WangJinXie2024}:
\begin{theorem}\label{thm:bounded_product_tnn}
Let $\Omega=\prod_{i=1}^{d}\Omega_i$ be a Cartesian-product domain, where each $\Omega_i\in \mathbb{R}$ is a bounded open interval. Consider a function $q\in H^m(\Omega)$. Assume that each subnetwork in the TNN is an FNN with one hidden layer and a non-polynomial activation function $\sigma\in C^m(\mathbb R)$. Then, for every $\delta>0$, there exist a rank $p\in\mathbb N$ and a TNN $\Psi_\theta$ such that
\begin{equation}\label{eq:bounded_product_tnn_approximation}
    \|q-\Psi_\theta\|_{H^m(\Omega)}<\delta.
\end{equation}   
\end{theorem}
\end{samepage}

We note that this theorem does not quantify the relation between separation rank $p$ and TNN the approximation error. Even for CP approximations, establishing a general rank--error relation remains an open research problem \cite{WangJinXie2024}.

This theorem does not directly apply to the distribution $f(t,x,v)$, since the velocity is defined over the unbounded domain $\mathbb R^3$. We therefore extend this approximation result to the domain $\mathcal D=[0,T]\times\Omega\times\mathbb{R}^3$ with $\Omega=\prod_{i=1}^{d}\Omega_i$ and $d\in\{1,2,3\}$. Let $\langle v\rangle=(1+|v|^2)^{1/2}$. For $s\geq0$ and $m\in\mathbb N_0$, we define a weighted Sobolev norm
\begin{equation}\label{eq:whole_velocity_weighted_norm}
    \|f\|_{H_s^m(\mathcal D)}^2
    :=\sum_{\alpha+|\beta|+|\gamma|\leq m}
    \left\|
        \langle v\rangle^s
        \partial_t^\alpha\partial_x^\beta\partial_v^\gamma f
    \right\|_{L^2(\mathcal D)}^2,
\end{equation}
where $\alpha\in\mathbb N_0$ is the time derivative order, and $\beta\in\mathbb N_0^d$ and $\gamma\in\mathbb N_0^3$ are multi-indices. We establish the TNN approximation theorem in this weighted Sobolev space:
\begin{samepage}
\begin{theorem}\label{thm:whole_velocity_weighted_tnn}
Let $f\in H_s^m(\mathcal D)$ be a particle density function. Then, for every $\delta>0$, there exists a rank $p\in\mathbb N$ and a TNN $f_\theta$ such that
\begin{equation}\label{eq:whole_velocity_tnn_approximation}
    \|f-f_\theta\|_{ H_s^m(\mathcal D)}<\delta.
\end{equation}
\end{theorem}
\end{samepage}

The proof can be found in Appendix~\ref{app:whole_velocity_tnn_proof}. For  $s>7/2$, the TNN $f_\theta$ can also approximate the velocity moments of $f$ up to second order, which will be needed to construct the loss. Specifically, the Cauchy--Schwarz inequality gives
\begin{equation}\label{eq:whole_velocity_moment_error}
\left\|
    \int_{\mathbb R^3}\varphi(v)(f-f_\theta)\,dv
\right\|_{L^2((0,T)\times\Omega)}
\leq C_s\|f-f_\theta\|_{H_s^m(\mathcal D)}
<C_s\delta,
\end{equation}
where $\varphi(v)\in\{1,v_i,v_iv_j,|v|^2\}$ with $i,j\in\{1,2,3\}$, and
$C_s=\bigl(\int_{\mathbb R^3}\langle v\rangle^{4-2s}\,dv\bigr)^{1/2}$ is finite when $s>7/2$.

In practice, each velocity factor in TNN $f_\theta$ is multiplied by a Gaussian envelope to enforce the tail decay of $f_\theta$ in velocity \cite{OhEtAl2025}
\begin{equation}\label{eq:gaussian_velocity_envelope}
    \widetilde{\xi}_r^{v_i}(v_i)
    =
    \xi_r^{v_i}(v_i)
    \exp\left[
        -\frac{\tau_i^2(v_i-\mu_i)^2}{2}
    \right],
    \quad
    i=1,\ldots,d_v,\quad r=1,\ldots,r_f.
\end{equation}
Here, $\mu_i$ and $\tau_i>0$ are trainable parameters initialized to $0$ and $1$, respectively. We still denote the modified velocity factors $\widetilde{\xi}_r^{v_i}$ by $\xi_r^{v_i}$ in what follows. Moreover, to improve numerical stability, each factor $\xi_r^s$ and $\eta_r^s$ in TNNs $f_\theta$ and $\phi_\theta$ is normalized in the $L^2$-norm, and the resulting scales are absorbed into $c_r^f$ and $c_r^\phi$, respectively.

The TNN architecture provides three computational advantages. First, the number of network forward passes is reduced from $\mathcal{O}(N^D)$ to $\mathcal{O}(DN)$. Second, the derivatives can be efficiently computed by forward-mode AD \cite{ChoEtAl2023}. Third, macroscopic velocity moments are reduced from high-dimensional integrals to sums of products of one-dimensional integrals. Specifically, for a multi-index $\boldsymbol{\alpha}=(\alpha_1,\ldots,\alpha_{d_v})\in\mathbb{N}_0^{d_v}$, we define $\mathbf{v}^{\boldsymbol{\alpha}}=\prod_{j=1}^{d_v}v_j^{\alpha_j}$ and evaluate the velocity moment as
\begin{equation}\label{eq:tnn_velocity_moment}
    \left\langle
        \mathbf{v}^{\boldsymbol{\alpha}}f_\theta
    \right\rangle
    =
    \sum_{r=1}^{r_f}
    c_r^f\,\xi_r^t(t)
    \prod_{i=1}^{d_x}\xi_r^{x_i}(x_i)
    \prod_{j=1}^{d_v}
    \left(
        \int_{\mathbb{R}}
        v_j^{\alpha_j}\xi_r^{v_j}(v_j)
        \,\mathrm{d}v_j
    \right).
\end{equation}
The computational cost of moments is reduced from $\mathcal{O}(N^{d_v})$ to $\mathcal{O}(d_vr_fN)$. This fast evaluation strategy will be used to construct the AP loss in Sect.~3.2.

\subsection{Flux closure-enhanced AP loss} 

We first review the concept of AP loss and explain how it captures macroscopic evolution in the limit. Following the PINN framework, the TNNs $f_\theta$ and $\phi_\theta$ can be trained  by minimizing a loss constructed by taking least square error of residuals of the kinetic system \eqref{eq:general_kinetic_model}. As $\varepsilon\to0$, the kinetic equations formally reduce to 
\begin{equation}\label{local-equil}
\mathbf{E}\cdot\nabla_{\mathbf{v}}f=\mathcal{Q}(f).
\end{equation}
This condition constrains only the local equilibrium and does not determine the macroscopic evolution for $\rho$ in the high-field limit. Hence the vanilla kinetic loss alone is insufficient to recover the correct high-field dynamics. To address this degeneracy in the limit, an asymptotic-preserving (AP) loss has been introduced for multiscale kinetic equations \cite{JinMaWu2023,LuWangXu2022}. Its key idea is to retain the original kinetic equation while incorporating macroscopic conservation laws into the loss. The resulting loss can thus preserve the essential macroscopic evolution as $\varepsilon\to0$. NNs trained with such losses are referred to as asymptotic-preserving neural networks (APNNs).

For \eqref{eq:general_kinetic_model}, the baseline AP loss can be constructed by incorporating the continuity equation \eqref{eq:macroscopic_mass_conservation} into the loss \cite{ZhangJin2025}. Accordingly, we need to introduce another TNN $\rho_\theta(t,\mathbf{x})$ to approximate the density $\rho$
\begin{equation}\label{eq:tnn_density}
    \rho_\theta(t,\mathbf{x})
    =
    \sum_{r=1}^{r_\rho}c_r^\rho
    \zeta_r^t(t)
    \prod_{i=1}^{d_x}\zeta_r^{x_i}(x_i),
\end{equation}
where $r_\rho$ is the separation rank and $\{c_r^\rho\}_{r=1}^{r_\rho}$ are trainable coefficients. The one-dimensional factors $\zeta_r^t$ and $\zeta_r^{x_i}$ are generated by independent subnetworks. The mass consistency condition $\rho_\theta=\langle f_\theta\rangle$ is also imposed in the loss. This AP loss provides a unified framework for the VPFP and BGK--Poisson systems. In particular, it does not require an explicit equilibrium distribution, which is not always available in the high-field regime.

We next examine how this baseline AP loss recovers the correct high-field limit. As $\varepsilon\to0$, the kinetic residual yields the local equilibrium \eqref{local-equil}. Taking the first velocity moment of this equation and using \eqref{eq:collision_first_moment}, we obtain the flux closure
\begin{equation}\label{eq:flux_closure}
     \mathbf{j}=\rho\mathbf{E}.
\end{equation}
Substituting this relation into the continuity equation \eqref{eq:macroscopic_mass_conservation}, together with the Poisson equation, recovers the high-field limit system \eqref{eq:high_field_limit}. The loss therefore has the AP property. However, the limiting flux closure \eqref{eq:flux_closure} is recovered from the kinetic residual without directly penalizing the closure error. This error can impair prediction accuracy for density evolution in the limit. 

To improve prediction accuracy in the high-field limit, we enhance the baseline AP loss by incorporating the first velocity-moment equation \eqref{eq:macroscopic_momentum_equation}
\begin{equation}\tag{\ref{eq:macroscopic_momentum_equation}}
    \varepsilon
    \left[
        \partial_t\langle\mathbf{v}f\rangle
        +\nabla_{\mathbf{x}}\cdot
        \langle\mathbf{v}\otimes\mathbf{v}f\rangle
    \right]
    +\langle\mathbf{v}f\rangle
    -\rho\mathbf{E}
    =0.
\end{equation}
As $\varepsilon\to0$, this equation reduces to the flux closure \eqref{eq:flux_closure}. The corresponding loss term therefore directly penalizes the closure error in the limit. Accordingly, we propose the following flux closure-enhanced AP formulation
\begin{equation}\label{eq:flux_closure_enhanced_ap_system}
\left\{\begin{aligned}
&\partial_t f +\mathbf{v}\cdot\nabla_{\mathbf{x}}f+\dfrac{1}{\varepsilon}\mathbf{E}\cdot\nabla_{\mathbf{v}}f=\dfrac{1}{\varepsilon}\mathcal{Q}(f), \\
&\partial_t\rho+\nabla_{\mathbf{x}}\cdot\langle\mathbf{v}f\rangle=0, \\
&\varepsilon\left[\partial_t\langle\mathbf{v}f\rangle+\nabla_{\mathbf{x}}\cdot\langle\mathbf{v}\otimes\mathbf{v}f\rangle\right]+\langle\mathbf{v}f\rangle-\rho\mathbf{E}=0, \\
&\rho=\langle f\rangle, \\
&-\Delta_{\mathbf{x}}\phi=\rho-h,\quad\mathbf{E}=-\nabla_{\mathbf{x}}\phi .\end{aligned}\right.
\end{equation}
The added equation significantly improves the prediction accuracy, as domestrated by the ablation study in Figure~\ref{fig:ap_errors}. The enhanced AP formulation also has three advantages. First, the first moment equation is obtained exactly from the kinetic equation for every $\varepsilon>0$. Its inclusion alters neither the kinetic equation nor its solutions. Second, since the required velocity moments are computed directly from $f$, no additional solution variables are introduced. Third, these moments can be efficiently evaluated using TNN structure via \eqref{eq:tnn_velocity_moment}. 

We now specify the loss function associated with system \eqref{eq:flux_closure_enhanced_ap_system}. Its three solution variables $f$, $\rho$, and $\phi$, are approximated by the TNNs $f_\theta(t,\mathbf{x},\mathbf{v})$, $\rho_\theta(t,\mathbf{x})$, and $\phi_\theta(t,\mathbf{x})$, respectively. Here, $\theta$ collects all trainable parameters and $\mathbf{E}_\theta=-\nabla_{\mathbf{x}}\phi_\theta$. The residuals of \eqref{eq:flux_closure_enhanced_ap_system} is defined as 
\begin{equation}\label{eq:aptnn_residuals}
\begin{aligned}
    \mathcal{R}_{\mathrm{kin}}
    &:={}
    \varepsilon\left(
        \partial_t f_{\theta}+\mathbf{v}\cdot\nabla_{\mathbf{x}}f_{\theta}\right)
    +\mathbf{E}_{\theta}\cdot\nabla_{\mathbf{v}}f_{\theta}-\mathcal{Q}_{\theta},\\
    \mathcal{R}_{\mathrm{cont}}&:={}\partial_t\rho_{\theta}+\nabla_{\mathbf{x}}\cdot\langle\mathbf{v}f_{\theta}\rangle,\\
    \mathcal{R}_{\mathrm{mom}}&:={}\varepsilon\left[\partial_t\langle\mathbf{v}f_{\theta}\rangle+\nabla_{\mathbf{x}}\cdot\langle\mathbf{v}\otimes\mathbf{v}f_{\theta}\rangle\right]
    +\langle\mathbf{v}f_{\theta}\rangle-\rho_{\theta}\mathbf{E}_{\theta},\\
    \mathcal{R}_{\mathrm{P}}&:={}-\Delta_{\mathbf{x}}\phi_{\theta}-\rho_{\theta}+h,\\
    \mathcal{R}_{\mathrm{mass}}&:={}\rho_{\theta}-\langle f_{\theta}\rangle,
\end{aligned}
\end{equation}
where $\mathcal{Q}_{\theta}=\mathcal{Q}_{\mathrm{FP}}(f_{\theta})$ for the VPFP system, whereas $\mathcal{Q}_{\theta}=\rho_{\theta}M-f_{\theta}$ for the BGK--Poisson system. Let $\Omega_{\mathbf{x}}=\prod_{i=1}^{d_x}\Omega_{x_i}$ denote the  bounded spatial domain and $\Omega_{\mathbf{v}}=\prod_{j=1}^{d_v}\Omega_{v_j}$ the truncated velocity domain for computation. The computational domains in both phase space and physical space are defined as 
\[
    \mathcal{D}_f
    =(0,T)\times\Omega_{\mathbf{x}}\times\Omega_{\mathbf{v}},
    \qquad
    \mathcal{D}_x
    =(0,T)\times\Omega_{\mathbf{x}},
\]
respectively. Their boundaries are denoted by $\Gamma_f=(0,T)\times\partial(\Omega_{\mathbf{x}}\times\Omega_{\mathbf{v}})$ and $\Gamma_\rho=\Gamma_\phi=(0,T)\times\partial\Omega_{\mathbf{x}}$. The losses for the equation residuals and IBCs are then defined using squared $L^2$-norms
\begin{equation}\label{eq:aptnn_loss_components}
\begin{aligned}
    \mathcal{L}_r(\theta)
    &:={}
    \lambda_{\mathrm{kin}}
        \|\mathcal{R}_{\mathrm{kin}}\|_{L^2(\mathcal{D}_f)}^2
    +\lambda_{\mathrm{cont}}
        \|\mathcal{R}_{\mathrm{cont}}\|_{L^2(\mathcal{D}_x)}^2
    +\lambda_{\mathrm{mom}}
        \|\mathcal{R}_{\mathrm{mom}}\|_{L^2(\mathcal{D}_x)}^2
    \\
    &\quad
    +\lambda_{\mathrm{P}}
        \|\mathcal{R}_{\mathrm{P}}\|_{L^2(\mathcal{D}_x)}^2
    +\lambda_{\mathrm{mass}}
        \|\mathcal{R}_{\mathrm{mass}}\|_{L^2(\mathcal{D}_x)}^2,
    \\
    \mathcal{L}_{\mathrm{bc}}(\theta)
    &:={}
    \|\mathcal{B}_f[f_{\theta}]\|_{L^2(\Gamma_f)}^2
    +\|\mathcal{B}_{\rho}[\rho_{\theta}]\|_{L^2(\Gamma_{\rho})}^2
    +\|\mathcal{B}_{\phi}[\phi_{\theta}]\|_{L^2(\Gamma_{\phi})}^2,
    \\
    \mathcal{L}_{\mathrm{ic}}(\theta)
    &:={}
    \|f_{\theta}(0)-f^0\|_{L^2(\Omega_{\mathbf{x}}
        \times\Omega_{\mathbf{v}})}^2
    +\|\rho_{\theta}(0)-\rho^0\|_{L^2(\Omega_{\mathbf{x}})}^2
    +\|\phi_{\theta}(0)-\phi^0\|_{L^2(\Omega_{\mathbf{x}})}^2.
\end{aligned}
\end{equation}
$\lambda_{\mathrm{kin}}$, $\lambda_{\mathrm{cont}}$, $\lambda_{\mathrm{mom}}$, $\lambda_{\mathrm{P}}$, and $\lambda_{\mathrm{mass}}$ are adjustable weights. The operators $\mathcal{B}_f$, $\mathcal{B}_{\rho}$, and $\mathcal{B}_{\phi}$ denote the BC operators for $f$, $\rho$ and $\phi$, respectively. The ICs are assumed to be compatible, with $\rho^0=\langle f^0\rangle$ and $-\Delta_{\mathbf{x}}\phi^0=\rho^0-h$. Finally, the total loss is defined as
\begin{equation}\label{eq:aptnn_loss}
    \mathcal{L}_{\mathrm{APTNN}}(\theta)
    =
    \mathcal{L}_r(\theta)
    +\lambda_{\mathrm{bc}}
        \mathcal{L}_{\mathrm{bc}}(\theta)
    +\lambda_{\mathrm{ic}}
        \mathcal{L}_{\mathrm{ic}}(\theta),
\end{equation}
where $\lambda_{\mathrm{bc}}$ and $\lambda_{\mathrm{ic}}$ are adjustable weights. We refer to the TNN-based method trained with this loss as the Asymptotic-Preserving Tensor Neural Network (APTNN). Due to the coordinate-separable structure of the TNN, the loss no longer needs to be evaluated on a full tensor-product grid. Instead, an efficient dimension-wise evaluation strategy will be presented in Sect.~\ref{sec:dimension_wise_loss}.

\subsection{Dimension-wise evaluation of the loss}\label{sec:dimension_wise_loss}
In this section, we show that the high-dimensional integrals in the APTNN loss \eqref{eq:aptnn_loss} can be decomposed into to a series of one-dimensional integrals. Directly approximating these integrals using the standard Gauss--Legendre (GL) quadrature rule incurs prohibitive computational costs and memory requirements. If $N$ quadrature points are used in each coordinate, evaluating the kinetic loss requires $N^{D_f}$ quadrature points. This curse of dimensionality is further exacerbated in UQ problems by additional random variables. To reduce this computational burden, we exploit the coordinate-separable structure of TNNs to evaluate the loss dimension-wise~\cite{WangJinXie2024,FengEtAl2025}.

We first show that the equation residuals \eqref{eq:aptnn_residuals} for loss construction inherit separable structure of TNNs. It is straightforward to verify that the operations used in these residuals, including products, velocity integration, partial differentiation and linear combinations of the resulting terms, preserve the separated structure. The Maxwellian $M$ in the BGK collision is also rank-one separable in velocity.  Therefore, $\mathcal{R}_{\mathrm{kin}}$, $\mathcal{R}_{\mathrm{cont}}$, each scalar component of $\mathcal{R}_{\mathrm{mom}}$, $\mathcal{R}_{\mathrm{P}}$ and $\mathcal{R}_{\mathrm{mass}}$ all admit separated representations. Let $\mathcal R(\mathbf y)$ denote any such residual component on the Cartesian-product domain $\mathcal D_{\mathcal R}=\prod_{\ell=1}^{D_{\mathcal R}}\Omega_\ell$, where $\mathbf y=(y_1,\ldots,y_{D_{\mathcal R}})$ and $\Omega_\ell$ is the one-dimensional domain of $y_\ell$. Here, $\mathcal D_{\mathcal R}$ represents $\mathcal D_f$ with dimension $D_{\mathcal R}=D_f$ for the kinetic residual, and $\mathcal D_x$ with dimension $D_{\mathcal R}=D_x$ for the remaining equation residuals. We can write $\mathcal R(\mathbf y)$ in a separated form:
\begin{equation}\label{eq:separated_residual}
    \mathcal R(\mathbf y)
    =\sum_{k=1}^{K_{\mathcal R}}a_k\prod_{\ell=1}^{D_{\mathcal R}}R_k^{(\ell)}(y_\ell),
\end{equation}
where $K_{\mathcal R}$ is the number of separated terms and $a_k$ are scalar coefficients.

We apply the separated form of $\mathcal R(\mathbf y)$ to decompose the equation loss $\mathcal L_r$ \eqref{eq:aptnn_loss_components} into a series of one-dimensional integrals. Specifically, the high-dimensional integral in each residual loss can be rewritten as:
\begin{equation}\label{eq:separated_residual_norm}
\begin{aligned}
    \mathcal L_{\mathcal R}
    &:=
    \|\mathcal R(\mathbf{y})\|_{L^2(\mathcal D_{\mathcal R})}^2
    =\int_{\prod_{\ell=1}^{D_{\mathcal R}}\Omega_\ell}
    \left(
        \sum_{k=1}^{K_{\mathcal R}}a_k
        \prod_{\ell=1}^{D_{\mathcal R}}
            R_k^{(\ell)}(y_\ell)
    \right)^2
    \,\mathrm dy_1\cdots\mathrm dy_{D_{\mathcal R}}\\
    &=\sum_{k=1}^{K_{\mathcal R}}
    \sum_{s=1}^{K_{\mathcal R}}a_ka_s\prod_{\ell=1}^{D_{\mathcal R}}
    \left[
        \int_{\Omega_\ell}
        R_k^{(\ell)}(y_\ell)
        R_s^{(\ell)}(y_\ell)
        \,\mathrm dy_\ell
    \right].
\end{aligned}
\end{equation}
One then approximates $\mathcal L_{\mathcal R}$ using the GL quadrature rule. For each $\ell=1,\ldots,D_{\mathcal R}$, let $\{y_n^{(\ell)},w_n^{(\ell)}\}_{n=1}^{N_\ell}$ denote the quadrature nodes and weights on $\Omega_\ell$. Each one-dimensional integral in \eqref{eq:separated_residual_norm} is approximated by
\begin{equation}\label{eq:residual_gram_matrix}
    G_{ks}^{(\ell)}
    :=
    \sum_{n=1}^{N_\ell}w_n^{(\ell)}
    R_k^{(\ell)}\!\left(y_n^{(\ell)}\right)
    R_s^{(\ell)}\!\left(y_n^{(\ell)}\right),
    \qquad
    k,s=1,\ldots,K_{\mathcal R}.
\end{equation}
Substituting \eqref{eq:residual_gram_matrix} into \eqref{eq:separated_residual_norm}, we obtain the dimension-wise approximation of $\mathcal L_{\mathcal R}$:
\begin{equation}\label{eq:dimension_wise_residual_loss}
    \mathcal L_{\mathcal R}
    \approx
    \sum_{k=1}^{K_{\mathcal R}}
    \sum_{s=1}^{K_{\mathcal R}}
    a_ka_s
    \prod_{\ell=1}^{D_{\mathcal R}}G_{ks}^{(\ell)}
    =\mathbf a^{\mathsf T}
    \left(
        G^{(1)}\odot\cdots\odot G^{(D_{\mathcal R})}
    \right)
    \mathbf a,
\end{equation}
where $G^{(\ell)}=[G_{ks}^{(\ell)}]_{k,s=1}^{K_{\mathcal R}} \in\mathbb R^{K_{\mathcal R}\times K_{\mathcal R}}$ is a weighted Gram matrix, $\mathbf a=(a_1,\ldots,a_{K_{\mathcal R}})^{\mathsf T}$ is the coefficient vector, and $\odot$ denotes the Hadamard product. Thus, all residual losses can be evaluated using one-dimensional quadrature and Gram-matrix contractions, and then assembled to form the complete equation loss $\mathcal L_r$. Here, separable BCs and ICs are considered, and their losses $\mathcal L_{\mathrm{bc}}$ and $\mathcal L_{\mathrm{ic}}$ are also evaluated dimension-wise. Finally, the total APTNN loss can be evaluated through one-dimensionals quadrature instead of direct tensor-product quadrature. 

We analyze the computational complexity and storage requirement for evaluating equation loss $\mathcal L_r$ dimension-wise using \eqref{eq:dimension_wise_residual_loss}. These estimates depend on the dimension $D_{\mathcal R}$, the number of separated terms $K_{\mathcal R}$ and the quadrature sizes $N_\ell$. Here, $K_{\mathcal R}$ is different for each equation residual in \eqref{eq:aptnn_residuals}, since it depends on the TNN ranks $r_f$, $r_\rho$ and $r_\phi$ together with the algebraic operations involved in each residual.  Let $K_{\mathrm{kin}}^{\mathrm{FP}}$, $K_{\mathrm{kin}}^{\mathrm{BGK}}$, $K_{\mathrm{cont}}$, $K_{\mathrm{mom},i}$, $K_{\mathrm{P}}$, and $K_{\mathrm{mass}}$ denote the numbers of separated terms in the corresponding equation residuals, where $K_{\mathrm{mom},i}$ corresponds to the $i$th scalar component of $\mathcal R_{\mathrm{mom}}$. These numbers satisfy the following upper bounds:
\begin{equation}\label{eq:residual_term_counts}
    \begin{aligned}
        K_{\mathrm{kin}}^{\mathrm{FP}}
        &\leq (1+d_x+d_v)r_f+d_v r_\phi r_f,
        &\qquad
        K_{\mathrm{kin}}^{\mathrm{BGK}}
        &\leq (d_x+2)r_f+d_v r_\phi r_f+r_\rho,\\
        K_{\mathrm{cont}}
        &\leq r_\rho+d_xr_f,
        &
        K_{\mathrm{mom},i}
        &\leq(d_x+2)r_f+r_\rho r_\phi,\\
        K_{\mathrm{P}}
        &\leq d_xr_\phi+r_\rho+1,
        &
        K_{\mathrm{mass}}
        &\leq r_\rho+r_f.
    \end{aligned}
\end{equation}
In particular, the products $E_{\theta,i}\partial_{v_i}f_\theta$ and $\rho_\theta E_{\theta,i}$ lead to multiplicative rank dependence into the separated term counts. They contribute at most $r_\phi r_f$ and $r_\rho r_\phi$ separated terms to $K_{\mathrm{kin}}$ and $K_{\mathrm{mom},i}$, respectively. The computational complexity and storage requirements for $\mathcal L_r$ is derived by combining the bounds \eqref{eq:residual_term_counts} with the loss evaluation formula \eqref{eq:dimension_wise_residual_loss}, as follows:

\begin{proposition}
\label{prop:dimension_wise_loss_complexity}
Assume that at most $N$ quadrature points are used in each coordinate. Let $K_{\mathrm{kin}}$, $K_{\mathrm{cont}}$, $K_{\mathrm{mom},i}$, $K_{\mathrm{P}}$ and $K_{\mathrm{mass}}$ satisfy the bounds in \eqref{eq:residual_term_counts}, where $K_{\mathrm{kin}}=K_{\mathrm{kin}}^{\mathrm{FP}}$ for the VPFP system and $K_{\mathrm{kin}}=K_{\mathrm{kin}}^{\mathrm{BGK}}$ for the BGK--Poisson system. When the loss of equation residual $\mathcal L_r$ in \eqref{eq:aptnn_loss_components} is evaluated dimension-wise using \eqref{eq:dimension_wise_residual_loss}, its computational complexity is
\begin{equation}\label{eq:aptnn_residual_loss_complexity}
    \mathcal O\!\left(
        N\left[
            D_f K_{\mathrm{kin}}^2
            +D_x\left(
                K_{\mathrm{cont}}^2
                +\sum_{i=1}^{d_v}K_{\mathrm{mom},i}^2
                +K_{\mathrm{P}}^2
                +K_{\mathrm{mass}}^2
            \right)
        \right]
    \right),
\end{equation}
where $D_f=1+d_x+d_v+d_z$ and $D_x=1+d_x+d_z$ with $d_z=0$ in the deterministic setting. The storage requirement for loss evaluation is 
\begin{equation}\label{eq:aptnn_residual_loss_storage}
\begin{aligned}
    \mathcal O\!\Bigg(
        &N\left[
            D_f K_{\mathrm{kin}}
            +D_x\left(
                K_{\mathrm{cont}}
                +\sum_{i=1}^{d_v}K_{\mathrm{mom},i}
                +K_{\mathrm{P}}
                +K_{\mathrm{mass}}
            \right)
        \right]\\
        &+K_{\mathrm{kin}}^2
        +K_{\mathrm{cont}}^2
        +\sum_{i=1}^{d_v}K_{\mathrm{mom},i}^2
        +K_{\mathrm{P}}^2
        +K_{\mathrm{mass}}^2
    \Bigg).
\end{aligned}
\end{equation}
\end{proposition}

\begin{proof}
We first estimate the computational cost. For any equation residual $\mathcal R(\mathbf y)$ of  the form \eqref{eq:separated_residual}, computing the Gram matrix $G^{(\ell)}$ in the $\ell$th coordinate used in \eqref{eq:dimension_wise_residual_loss} requires $\mathcal O(N_\ell K_{\mathcal R}^2)$ operations. Since $N_\ell\leq N$, constructing all $D_{\mathcal R}$ Gram matrices requires $\mathcal O(ND_{\mathcal R}K_{\mathcal R}^2)$ operations. The Hadamard product and the final quadratic form require $\mathcal O(D_{\mathcal R}K_{\mathcal R}^2)$ and $\mathcal O(K_{\mathcal R}^2)$ operations, respectively. Since these two costs are of no higher order than the Gram-matrix construction cost, summing the bounds $\mathcal O(ND_{\mathcal R}K_{\mathcal R}^2)$ over all residual components in $\mathcal L_r$ yields the result.

We next estimate the storage required for loss evaluation. For any residual $\mathcal R(\mathbf y)$, storing the values $R_k^{(\ell)}(y_n^{(\ell)})$ for Gram-matrix computation in \eqref{eq:residual_gram_matrix} requires $N_\ell K_{\mathcal R}$ entries per coordinate, giving a total storage of $\mathcal O(ND_{\mathcal R}K_{\mathcal R})$. The Hadamard product in \eqref{eq:dimension_wise_residual_loss} can be accumulated sequentially across coordinates, retaining only the current Gram matrix and the accumulated product matrix. Both matrices have $K_{\mathcal R}^2$ entries, so the additional storage is $\mathcal O(K_{\mathcal R}^2)$. Summing these storage bounds over all residual components in $\mathcal L_r$ yields \eqref{eq:aptnn_residual_loss_storage}.
\end{proof}

To extract the leading-order rank dependence of both bounds in Proposition~\ref{prop:dimension_wise_loss_complexity}, we assume that the TNN ranks satisfy $\min\{r_f,r_\rho,r_\phi\}\gg\max\{1,d_x,d_v\}$. Since $\rho$ and $\phi$ are independent of velocity, we choose smaller ranks $r_\rho,r_\phi<r_f$. Under the above assumptions, both bounds are controlled by the kinetic residual contribution with $K_{\mathrm{kin}}=\mathcal O(d_vr_fr_\phi)$. Hence, these bounds simplify as follows:
\begin{corollary}\label{cor:aptnn_rank_complexity}
Under the rank assumptions above, the computational complexity of the dimension-wise evaluation of the equation loss $\mathcal L_r$ is
\begin{equation}\label{eq:aptnn_rank_complexity}
    \mathcal O\!\left(
        ND_f d_v^2r_f^2r_\phi^2
    \right).
\end{equation}
The corresponding storage requirement is 
\begin{equation}\label{eq:aptnn_rank_storage}
    \mathcal O\!\left(
        ND_f d_vr_fr_\phi
        +d_v^2r_f^2r_\phi^2
    \right).
\end{equation}
For fixed TNN ranks and quadrature size, the computational and storage bounds have at most cubic and quadratic dependence on the physical dimensions, respectively, and at most linear dependence on the random dimension.
\end{corollary}

In contrast, for previous TNN methods applied to kinetic equations~\cite{OhEtAl2025,ChenDimarcoPareschi2025}, full-grid evaluation of the equation loss entails $\mathcal O(N^{D_f})$ computational complexity and storage requirements.

\section{Numerical results}\label{sec:numerical_results}

In this section, we present numerical experiments to assess the accuracy, computational efficiency and scalability of APTNN. Reference solutions are computed using the AP schemes developed for the VPFP system \cite{JinWang2011} and the semiconductor Boltzmann equation \cite{JinWang2013}. We extend both schemes to the four- and six-dimensional phase spaces and implement them in Python. Both schemes employ a second-order IMEX time discretization, a MUSCL spatial discretization, and an FFT-based Poisson solver. 

In APTNN numerical realization, each one-dimensional subnetwork adopts the modified MLP architecture \cite{WangTengPerdikaris2021} with the Swish activation function \cite{RamachandranEtAl2017}. The subnetworks for the physical coordinates have $4$ hidden layers of $64$ neurons each, and those for the random coordinates have $3$ hidden layers of $32$ neurons. All weight matrices are initialized using the Glorot normal initializer. Periodic BCs are used in each spatial direction for all examples. These conditions are enforced as hard constraints by applying a Fourier embedding to each spatial input \cite{LuEtAl2022,LyuEtAl2021}. Specifically, each input \(x_i\in[a_i,b_i]\) is mapped to \(\gamma_i(x_i)=\bigl(\cos(\omega_i x_i),\sin(\omega_i x_i),\ldots,\cos(J\omega_i x_i),\sin(J\omega_i x_i)\bigr)\), where \(\omega_i=2\pi/P_i\) and \(P_i=b_i-a_i\). We set \(J=1\) throughout. Accordingly, the boundary loss \(\mathcal{L}_{\mathrm{bc}}\) in \eqref{eq:aptnn_loss} is omitted. The remaining loss terms are evaluated with $N$ Gauss--Legendre (GL) quadrature points in each coordinate direction. The networks are trained for $60{,}000$ iterations using the Adam optimizer. The learning rate is initialized to $10^{-3}$ and reduced by a factor of $0.9$ every $1{,}000$ iterations. The implementation uses JAX and Flax.

We report relative space--time $L^2$-errors for the density $\rho$ and the electric field $\mathbf{E}$. In particular, we focus on density perturbations around the uniform background $h=1$ and define the relative error in $\rho-1$ as
\begin{equation}
    e_{\rho-1}
    =
    \frac{
        \left\|\rho_{\mathrm{pred}}-\rho_{\mathrm{ref}}\right\|_{L^2([0,T]\times\Omega_{\mathbf{x}})}
    }{
        \left\|\rho_{\mathrm{ref}}-1\right\|_{L^2([0,T]\times\Omega_{\mathbf{x}})}
    },
    \label{eq:rho_relative_error}
\end{equation}
where $\rho_{\mathrm{pred}}$ and $\rho_{\mathrm{ref}}$ denote the predicted and reference densities, respectively. The error $e_{\mathbf{E}}$ of the electric field is defined analogously. For multi-dimensional problems, the norm of $\mathbf{E}=(E_1,\ldots,E_{d_x})$ is defined by $\|\mathbf{E}\|_{L^2}^2=\sum_{i=1}^{d_x}\|E_i\|_{L^2}^2$. The reference solutions are generated on a single node of an HPC cluster with two Intel Xeon Platinum 8358 processors (32 cores each, 2.6\,GHz) and 512\,GB of memory. All network training is performed on a single NVIDIA A100 GPU with 40\,GB of memory. We note that no extensive hyperparameter tuning is performed for the networks.

\subsection{AP property and prediction accuracy} 

We consider a one-dimensional Landau damping problem to examine the AP property of APTNN for both the VPFP system and the BGK--Poisson system. The initial distribution is given by
\begin{equation}\label{eq:landau_ic}
f^0(x,v)=\big(1+\alpha\cos(kx)\big)\frac{e^{-v^2/2}}{\sqrt{2\pi}},
\end{equation}
where $\alpha=0.1$ and $k=0.25$. The spatial domain is $x\in[0,L]$ with $L=2\pi/k$ and the velocity domain is truncated to $v\in[-V,V]$ with $V=6$. We set the terminal time $T=1$ and consider $\varepsilon\in\{1,10^{-1},10^{-2},10^{-3}\}$. Reference solutions for both models are computed using $256$ uniform spatial cells, $128$ velocity grid points and a time step $\Delta t=2\times10^{-3}$. We compare vanilla TNN, baseline APTNN and the proposed APTNN. The vanilla TNN loss is based on the original kinetic system, whereas the proposed APTNN adds the first velocity-moment equation to the baseline AP loss. We set the separation ranks to $r_f=50$ and $r_\phi=20$ for all three methods, and $r_\rho=20$ for both APTNN variants. The loss is evaluated with $N=32$ GL quadrature points in each coordinate direction, and all loss weights are set to one. With these settings, training the proposed APTNN takes approximately $4$ minutes per run for both the VPFP and BGK--Poisson systems.

\begin{figure}[htbp]
    \centering
    \begin{minipage}[t]{\textwidth}
        \centering
        \includegraphics[width=0.65\linewidth]
            {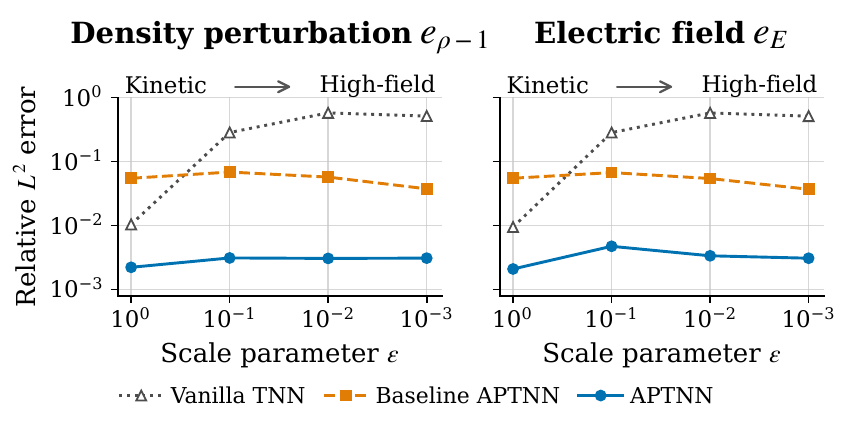}
        \par
        {\footnotesize (a) VPFP system}
    \end{minipage}
    \par\medskip
    \begin{minipage}[t]{\textwidth}
        \centering
        \includegraphics[width=0.65\linewidth]
            {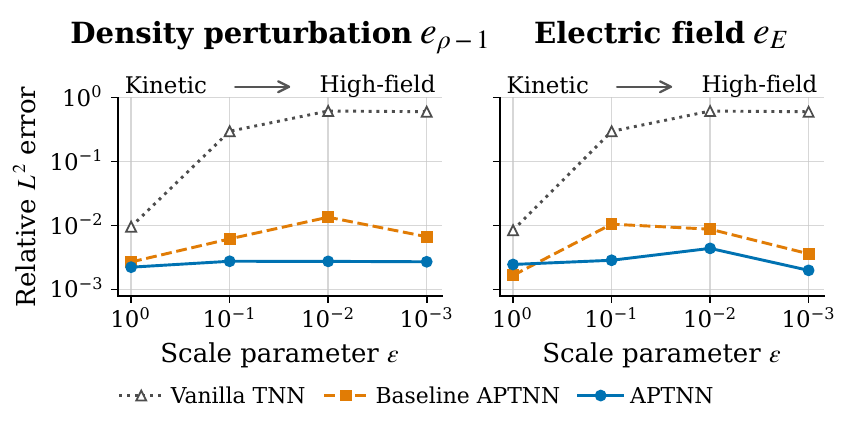}
        \par
        {\footnotesize (b) BGK-Poisson system}
    \end{minipage}
    \caption{AP property and prediction accuracy for the one-dimensional Landau damping problem. Relative space--time $L^2$-errors of the density perturbation $\rho-1$ and the electric field $E$ are shown for (a) the VPFP system and (b) the BGK--Poisson system. The three methods are vanilla TNN, baseline APTNN and the proposed APTNN, which adds the first velocity-moment equation to the baseline AP loss. The scale parameter $\varepsilon$ decreases from left to right, from the kinetic to the high-field regime.}
    \label{fig:ap_errors}
\end{figure}

Figure~\ref{fig:ap_errors} compares the relative $L^2$-errors of the density perturbation $\rho-1$ and the electric field $E$ obtained using the three loss formulations. The vanilla TNN becomes inaccurate as $\varepsilon$ decreases and therefore does not have the AP property. In contrast, both APTNN variants remain accurate for all tested $\varepsilon$, supporting their AP property. Incorporating the first velocity-moment equation into the AP loss further improves prediction accuracy. Specifically, for the VPFP system, the proposed APTNN reduces the errors from $\mathcal{O}(10^{-2})$ to $\mathcal{O}(10^{-3})$. For the BGK--Poisson system, the baseline APTNN errors fluctuate between $10^{-3}$ and $10^{-2}$,  whereas the proposed APTNN stays around $10^{-3}$. 

\begin{figure}[htbp]
    \centering
    \begin{minipage}[t]{0.48\textwidth}
        \centering
        \includegraphics[width=\linewidth]
            {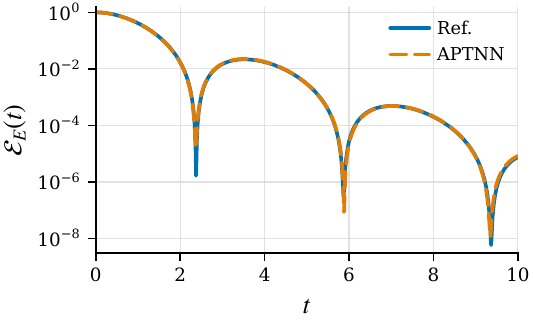}
        \par
        {\footnotesize (a) VPFP system}
    \end{minipage}\hspace{0.02\textwidth}
    \begin{minipage}[t]{0.48\textwidth}
        \centering
        \includegraphics[width=\linewidth]
            {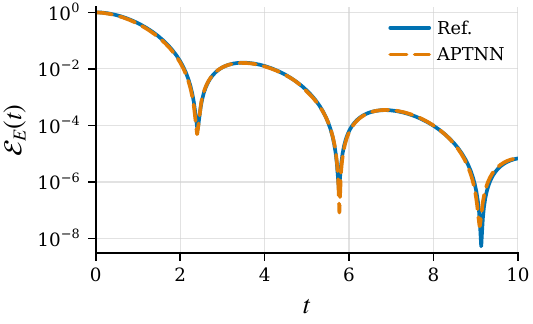}
        \par
        {\footnotesize (b) BGK--Poisson system}
    \end{minipage}
    \caption{Evolution of the electric-field energy $\mathcal{E}_E(t)$ for the one-dimensional Landau damping problem with $T=10$ and $\varepsilon=1$ for (a) the VPFP system and (b) the BGK--Poisson system. The solid and dashed curves denote the reference solutions and the APTNN predictions, respectively. The vertical axis is on a logarithmic scale.}
    \label{fig:electric_field_energy}
\end{figure}

We next extend the simulations to $T=10$ with $\varepsilon=1$ to examine the evolution of the electric-field energy, defined by
\begin{equation}
    \mathcal{E}_E(t)
    =
    \frac{1}{2}
    \int_{\Omega_x}|E(t,x)|^2\,\mathrm{d}x
    =
    \frac{1}{2}
    \|E(t,\cdot)\|_{L^2(\Omega_x)}^2.
    \label{eq:electric_field_energy}
\end{equation}
To stabilize training, we divide $[0,10]$ into two windows, $[0,5]$ and $[5,10]$, and train sequentially using a time-marching strategy~\cite{ChoEtAl2023}. The terminal prediction from the first window is used as the ICs for the second. The network ranks and quadrature settings remain unchanged. The loss weights for the kinetic, continuity and Poisson residuals are adjusted to $50$ for both systems, while all other weights remain one. The network is trained for $60{,}000$ iterations on each window, with a training time comparable to that of the $T=1$ case. Figure~\ref{fig:electric_field_energy} shows good agreement between the APTNN predictions and the reference solutions for both systems. The relative $L^2([0,T])$-errors of the electric-field energy are $4.48\times10^{-3}$ for the VPFP system and $3.06\times10^{-2}$ for the BGK-Poisson system. 

\subsection{Efficiency of dimension-wise loss evaluation}\label{sec:2d2v_gpu} 

This section assess the computational cost and memory usage of dimension-wise loss evaluation. It also checks the ability of APTNN to capture non-equilibrium dynamics. We consider a 2D2V VPFP problem with a double-stream initial distribution
\begin{equation}\label{eq:2d2v_double_stream_ic}
    f^0(\mathbf{x},\mathbf{v})
    =
    \left[1+\alpha\cos(kx_1)\cos(kx_2)\right]
    \frac{
        e^{-(v_1-1.5)^2/2}+e^{-(v_1+1.5)^2/2}
    }{2\sqrt{2\pi}}
    \frac{e^{-v_2^2/2}}{\sqrt{2\pi}},
\end{equation}
where $\alpha=0.05$ and $k=0.5$. The spatial domain is $(x_1,x_2)\in[0,2\pi/k]^2$ and the velocity domain is truncated to $(v_1,v_2)\in[-6,6]^2$. We set $T=1$ and $\varepsilon=10^{-3}$. The reference solution is computed using $128$ uniform cells in each spatial direction, $48$ grid points in each velocity direction and a time step $\Delta t=2\times10^{-3}$.

We compare dimension-wise and full-grid evaluation of the same loss. The former uses one-dimensional quadrature and Gram-matrix contractions, whereas the latter evaluates the loss directly using full tensor-product quadrature points. The quadrature orders are set as $N=8,12,16,24,32,40,48,64$. The separation ranks of TNNs are set to $r_f=50$ and $r_\rho=r_\phi=20$, and all loss weights to one. For each $N$, both strategies are evaluated using the same network parameters. We time $30$ synchronized training steps after JIT compilation and warm-up, and report the median wall-clock time per step and peak GPU memory usage. For all tested $N$, the absolute difference between the loss values computed by the two strategies is of order $10^{-13}$. Hence, we report only the relative $L^2$-errors of APTNN trained by dimension-wise strategy.

\begin{figure}[htbp]
    \centering
    \begin{minipage}[t]{0.325\textwidth}
        \centering
        \includegraphics[width=\linewidth]
            {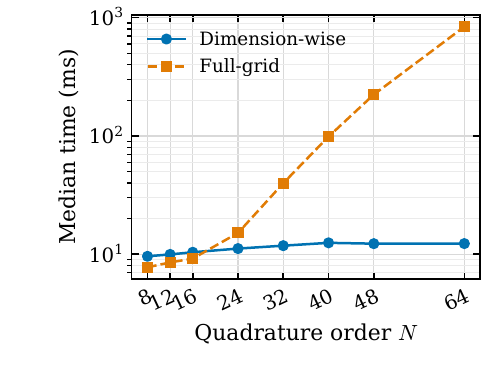}
        \par\vspace{-0.3em}
        {\footnotesize (a) Time per training step}
    \end{minipage}
    \hfill
    \begin{minipage}[t]{0.325\textwidth}
        \centering
        \includegraphics[width=\linewidth]
            {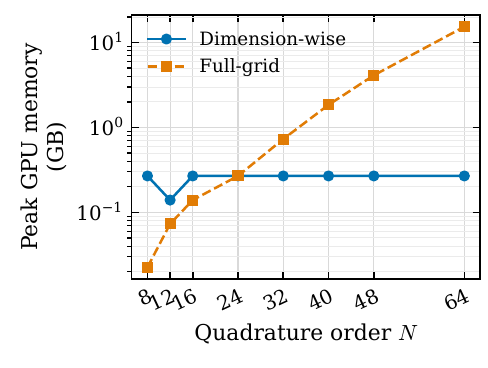}
        \par
        {\footnotesize (b) Peak GPU memory}
    \end{minipage}
    \hfill
    \begin{minipage}[t]{0.325\textwidth}
        \centering
        \includegraphics[width=\linewidth]
            {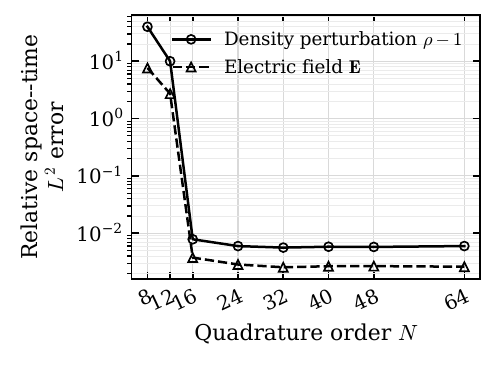}
        \par
        {\footnotesize (c) APTNN relative errors}
    \end{minipage}
    \caption{Computational cost, memory usage and errors for the 2D2V double-stream problem with $\varepsilon=10^{-3}$ and $T=1$ as the quadrature order $N$ varies. Panels~(a) and~(b) show the median training-step time and peak GPU memory, respectively, for both dimension-wise and full-grid loss evaluation. Panel~(c) shows the relative space--time $L^2$ errors of the density perturbation $\rho-1$ and the electric field $\mathbf{E}$ for APTNN trained by dimension-wise evaluation.}
    \label{fig:2d2v_training_scaling}
\end{figure}

Figure~\ref{fig:2d2v_training_scaling} exhibits the time per training step, peak GPU memory usage and relative $L^2$-errors in panels~(a)--(c), respectively. Panels~(a) and (b) show that dimension-wise evaluation maintains relatively stable time and memory costs over the tested range, requiring at most $12.5$\,ms per training step and $0.27$\,GB of peak GPU memory. In contrast, full-grid evaluation exhibits a steep increase in both metrics as $N$ increases. Since this strategy computs the kinetic loss using all $N^5$ quadrature points, its time per training step rises from $7.77$\,ms to $838.13$\,ms, and its peak GPU memory usage rises from $0.023$\,GB to $15.37$\,GB. For $N\le16$, the full-grid evaluation is slightly more efficient because the dimension-wise evaluation is dominated by the fixed rank-dependent cost of Gram-matrix contractions. However, the errors have not yet converged in this range as shown in panel~(c). From $N=24$ onward, dimension-wise evaluation becomes more efficient and its advantage grows rapidly with $N$. 

These results demonstrate that dimension-wise evaluation requires substantially less computation and memory than full-grid evaluation to support fine quadrature for accurate predictions. This advantage becomes more pronounced in higher-dimensional kinetic and UQ problems, because the $N^{D_f}$ full-grid size makes increasing $N$ much more expensive.

\begin{figure}[htbp]
    \centering
    \includegraphics[width=\textwidth]
        {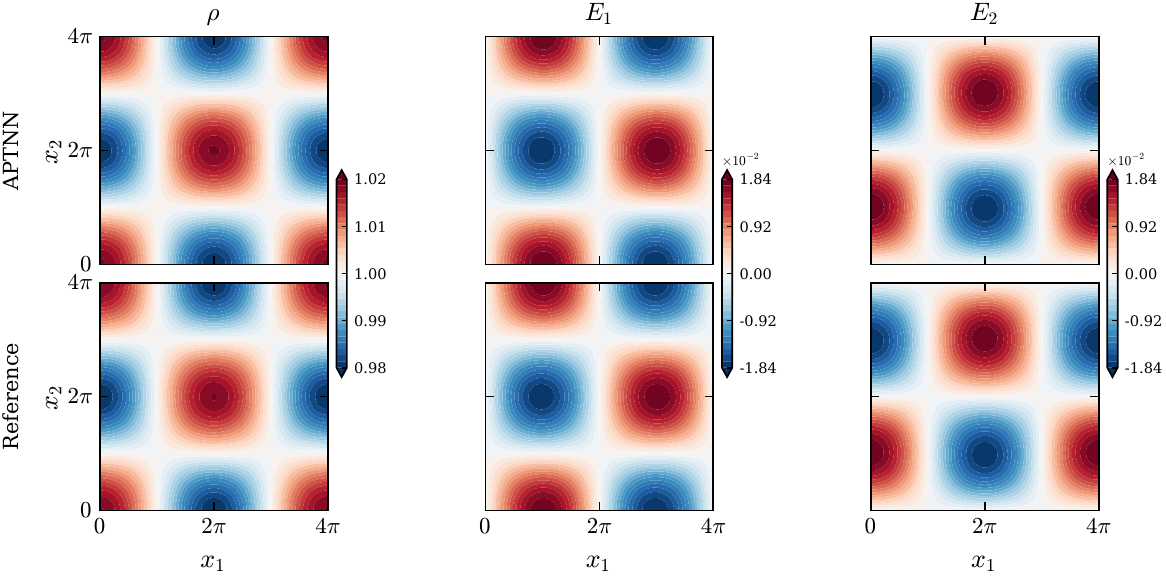}
    \caption{Results for the 2D2V double-stream problem with $\varepsilon=10^{-3}$ at $T=1$ using $N=32$. The top and bottom rows show the APTNN predictions and reference solutions, respectively. The columns show the density $\rho$ and the electric-field components $E_1$ and $E_2$, respectively.}
    \label{fig:2d2v_fields_T1}
\end{figure}

Figure~\ref{fig:2d2v_fields_T1} compares the APTNN predictions for $\rho$, $E_1$ and $E_2$ with the reference solution at $T=1$ using $N=32$. The relative space--time $L^2$-errors of the density perturbation $\rho-1$ and the vector electric field $\mathbf{E}=(E_1,E_2)$ are $5.62\times10^{-3}$ and $2.55\times10^{-3}$, respectively. Therefore, APTNN can capture the non-equilibrium dynamics accurately. The training time is approximately $7$ minutes. It is measured without per-step synchronization and is not directly comparable to the time per training step in Figure~\ref{fig:2d2v_training_scaling}(a).

\subsection{3D3V problem}
We examine the accuracy of APTNN in high-dimensional settings using a 3D3V BGK--Poisson problem. The initial distribution is
\begin{equation}\label{eq:3d3v_smooth_ic}
    f^0(\mathbf{x},\mathbf{v})
    =
    \left[1+\alpha\sum_{i=1}^{3}\cos(kx_i)\right]
    \prod_{i=1}^{3}\frac{e^{-v_i^2/2}}{\sqrt{2\pi}},
\end{equation}
where $\alpha=0.05$ and $k=0.5$. The spatial domain is $(x_1,x_2,x_3)\in[0,2\pi/k]^3$ and the velocity domain is truncated to $(v_1,v_2,v_3)\in[-6,6]^3$. We set $T=1$ and consider $\varepsilon\in\{1,10^{-3}\}$. The reference solutions are computed using $80$ uniform cells in each spatial direction, $16$ grid points in each velocity direction and a time step $\Delta t=2\times10^{-3}$. The separation ranks of TNN are set to $r_f=50$ and $r_\rho=r_\phi=20$. The loss is evaluated using $N=32$ GL quadrature points in each coordinate direction. For $\varepsilon=1$, the IC loss for $f$ is assigned a weight of $200$, while all other loss weights are set to $1$. For $\varepsilon=10^{-3}$, all loss weights are set to $1$. The corresponding training times are approximately $12$ and $14$ minutes, respectively.

\begin{figure}[htbp]
    \centering
    \includegraphics[width=\textwidth]
        {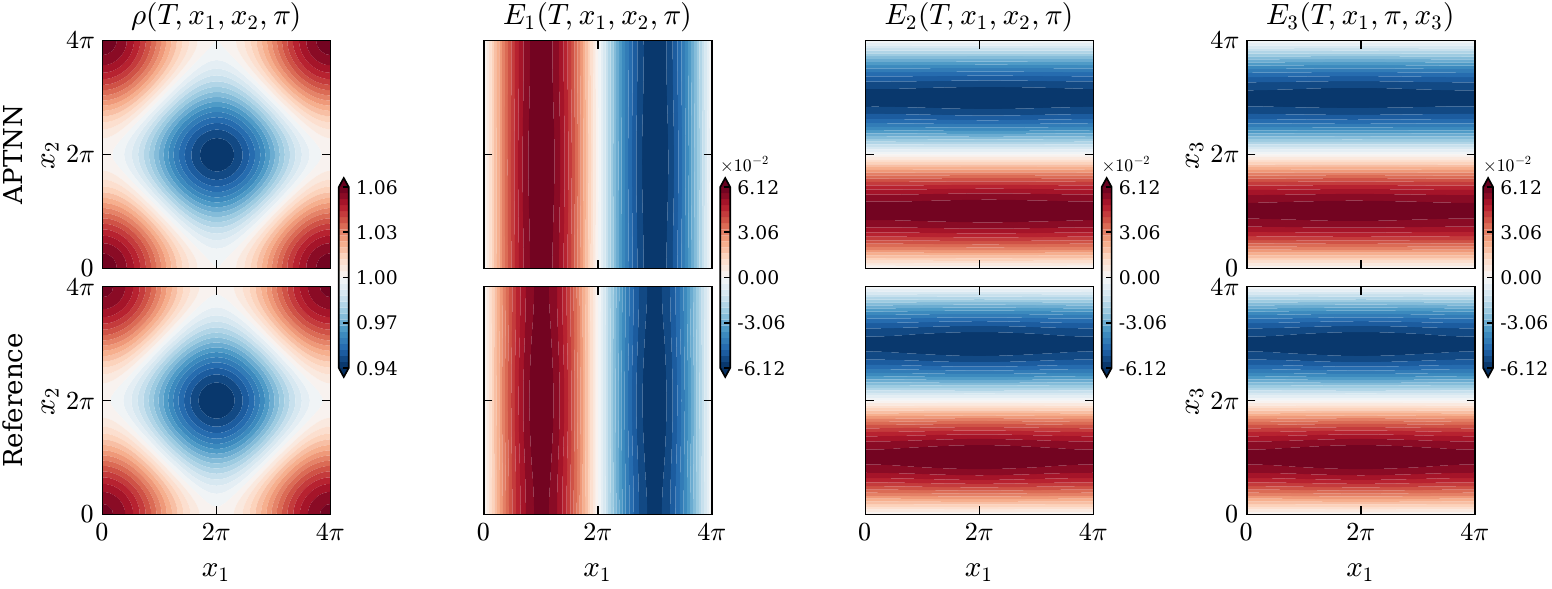}
    \caption{Results for the 3D3V BGK--Poisson problem with $\varepsilon=1$ at $T=1$. The top and bottom rows show the APTNN predictions and reference solutions, respectively. The columns for $\rho$, $E_1$, and $E_2$ show the $x_1$--$x_2$ plane at $x_3=\pi$, while the $E_3$ column shows the $x_1$--$x_3$ plane at $x_2=\pi$.}
    \label{fig:3d3v_bgk_fields_eps1_T1}
\end{figure}

\begin{figure}[htbp]
    \centering
    \includegraphics[width=\textwidth]
        {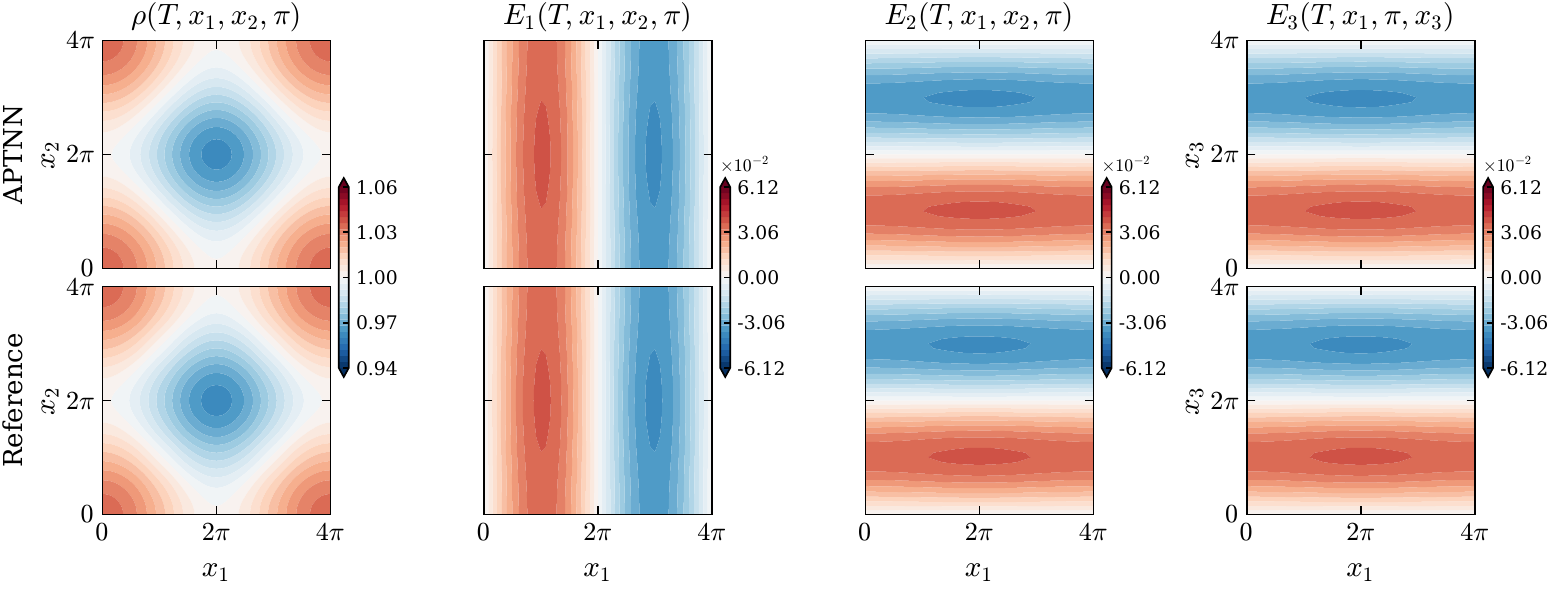}
    \caption{Results for the 3D3V BGK--Poisson problem with $\varepsilon=10^{-3}$ at $T=1$. The top and bottom rows show the APTNN predictions and reference solutions, respectively. The columns for $\rho$, $E_1$, and $E_2$ show the $x_1$--$x_2$ plane at $x_3=\pi$, while the $E_3$ column shows the $x_1$--$x_3$ plane at $x_2=\pi$.}
    \label{fig:3d3v_bgk_fields_eps1e-3_T1}
\end{figure}

Figures~\ref{fig:3d3v_bgk_fields_eps1_T1} and \ref{fig:3d3v_bgk_fields_eps1e-3_T1} compare the APTNN predictions with the reference solutions at $T=1$ for $\varepsilon=1$ and $10^{-3}$, respectively. The predicted density and electric field agree well with the reference solutions in both regimes. For $\varepsilon=1$, the relative space--time $L^2$-errors of the density perturbation $\rho-1$ and the vector electric field $\mathbf{E}=(E_1,E_2,E_3)$ are $8.78\times10^{-3}$ and $7.02\times10^{-3}$, respectively. The corresponding errors for $\varepsilon=10^{-3}$ are $1.98\times10^{-3}$ and $1.61\times10^{-3}$. 

\subsection{UQ for 3D3V problems}

We consider the 3D3V VPFP system with random inputs to evaluate the performance of APTNN in multi-dimensional UQ problems. Two cases with random dimensions $d_z=3$ and $d_z=5$ are tested. For $d_z=3$, the initial distribution is 
\begin{equation}\label{eq:3d3v_uq_d3_ic}
    f^0_{3}(\mathbf{x},\mathbf{v},\mathbf{z})
    =
    \left[
        1+\alpha\sum_{i=1}^{3}
        (1+z_i)\cos(kx_i)
    \right]
    \prod_{i=1}^{3}\frac{e^{-v_i^2/2}}{\sqrt{2\pi}}.
\end{equation}
For $d_z=5$, we take 
\begin{equation}\label{eq:3d3v_uq_d5_ic}
    f^0_{5}(\mathbf{x},\mathbf{v},\mathbf{z})
    =
    \left[
        1+\alpha
        \left(1+\frac{1}{5}\sum_{m=1}^{5}z_m\right)
        \sum_{i=1}^{3}\cos(kx_i)
    \right]
    \prod_{i=1}^{3}\frac{e^{-v_i^2/2}}{\sqrt{2\pi}}.
\end{equation}
Here, $z_1,\ldots,z_{d_z}$ are independent random variables uniformly distributed on $[-1,1]$. The spatial and velocity domains, BCs and parameters $\alpha$ and $k$ remain the same as in the preceding 3D3V problem. We set $\varepsilon=10^{-4}$ and $T=1$ and evaluate the means and standard deviations of $\rho$ and $\mathbf{E}$. The reference statistics are computed nonintrusively using a level-4 Smolyak sparse-grid rule based on nested Clenshaw--Curtis points \cite{Smolyak1963,ClenshawCurtis1960}. The deterministic solver uses $48$ uniform cells in each spatial direction, $20$ grid points in each velocity direction and a time step $\Delta t=4\times10^{-3}$. The separation ranks of TNNs are set to $r_f=50$ and $r_\rho=r_\phi=30$. The loss is evaluated using $32$ GL quadrature points in each physical coordinate and $16$ in each random coordinate. All loss weights are set to $1$. The training time is approximately $23$ minutes for $d_z=3$ and $25$ minutes for $d_z=5$. In contrast, full-grid loss evaluation would far exceed the memory capacity of a single GPU in both cases. Even for $d_z=3$, evaluating the kinetic loss would require over $500$\,TB of memory.

\begin{figure}[htbp]
    \centering
    \includegraphics[width=\textwidth]
        {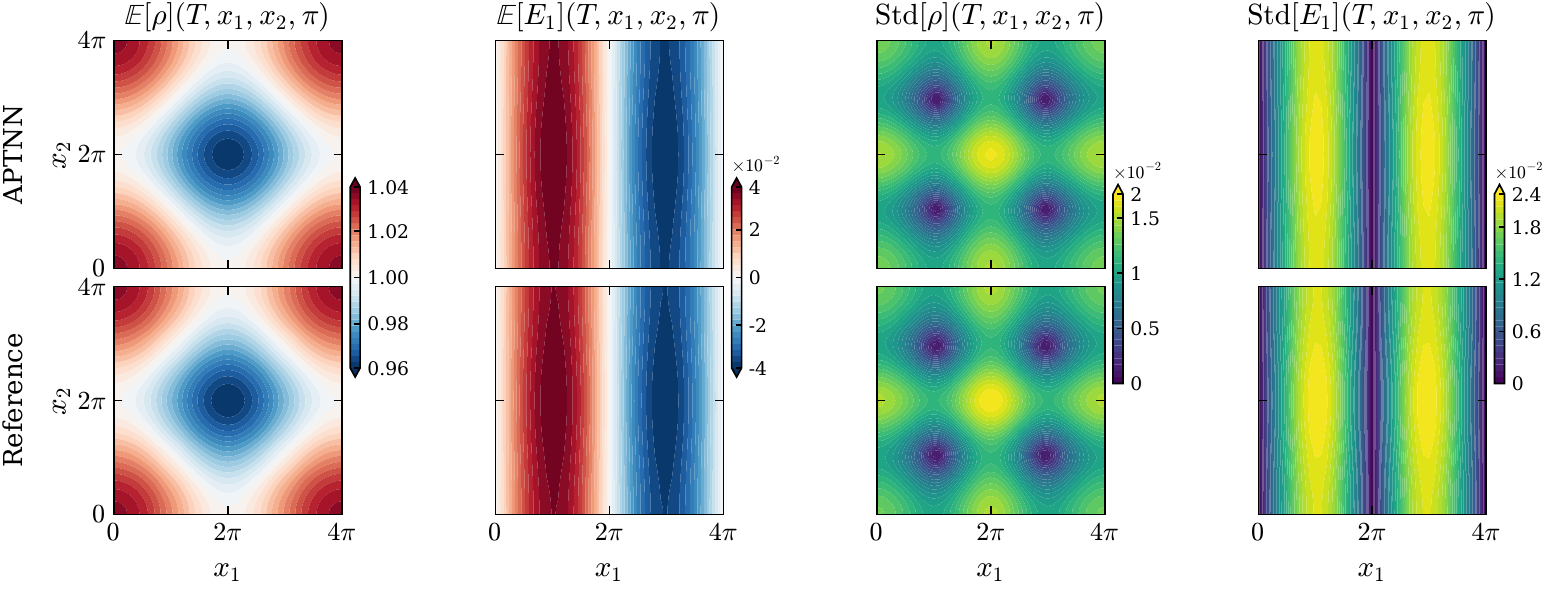}
    \caption{Means and standard deviations in the UQ problem for the 3D3V VPFP system with $d_z=3$ and $\varepsilon=10^{-4}$ at $T=1$. The top and bottom rows show the APTNN predictions and reference statistics, respectively. The columns show $\mathbb{E}[\rho]$, $\mathbb{E}[E_1]$, $\operatorname{Std}[\rho]$ and $\operatorname{Std}[E_1]$ on the $x_1$--$x_2$ plane at $x_3=\pi$.}
    \label{fig:3d3v_vpfp_uq_dz3_fields}
\end{figure}

\begin{figure}[htbp]
    \centering
    \includegraphics[width=\textwidth]
        {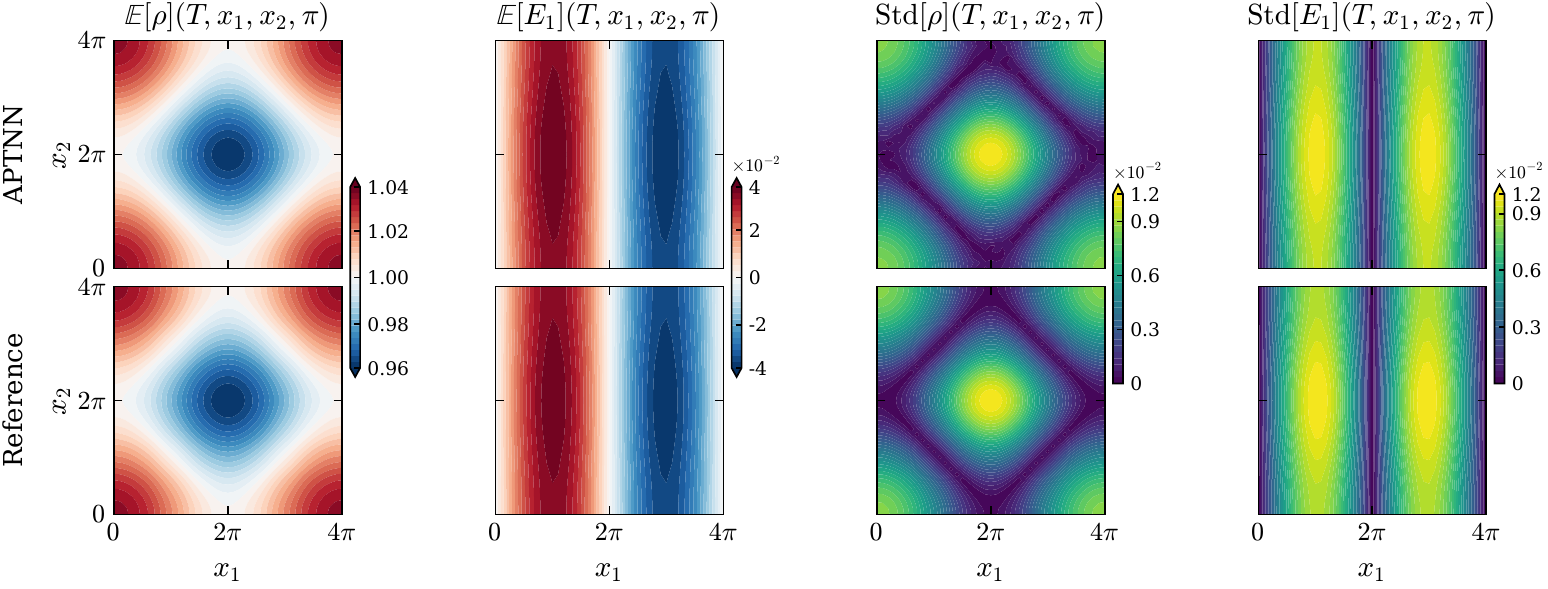}
    \caption{Means and standard deviations in the UQ problem for the 3D3V VPFP system with $d_z=5$ and $\varepsilon=10^{-4}$ at $T=1$. The top and bottom rows show the APTNN predictions and reference statistics, respectively. The columns show $\mathbb{E}[\rho]$, $\mathbb{E}[E_1]$, $\operatorname{Std}[\rho]$ and $\operatorname{Std}[E_1]$ on the $x_1$--$x_2$ plane at $x_3=\pi$.}
    \label{fig:3d3v_vpfp_uq_dz5_fields}
\end{figure}

Figures~\ref{fig:3d3v_vpfp_uq_dz3_fields} and \ref{fig:3d3v_vpfp_uq_dz5_fields} compare the APTNN predictions with the reference statistics at $T=1$ for $d_z=3$ and $d_z=5$, respectively. The predicted means and standard deviations of both the density and the electric field agree closely with the reference statistics in both cases. For $d_z=3$, the relative space--time $L^2$-errors for the density are $3.87\times10^{-3}$ for $\mathbb{E}[\rho]-1$ and $1.00\times10^{-2}$ for $\operatorname{Std}[\rho]$, while those for the electric field are $2.51\times10^{-3}$ for $\mathbb{E}[\mathbf{E}]$ and $7.31\times10^{-3}$ for $\operatorname{Std}[\mathbf{E}]$. For $d_z=5$, the corresponding errors are $3.81\times10^{-3}$ and $1.47\times10^{-2}$ for the density, and $2.44\times10^{-3}$ and $1.13\times10^{-2}$ for the electric field.

Finally, Table~\ref{tab:training_times} summarizes the wall-clock training times for the representative experiments presented in Sect.~\ref{sec:numerical_results}.

\begin{table}[htbp]
    \centering
    \footnotesize
    \caption{
        Wall-clock training times of APTNN for the numerical experiments in Sect.~\ref{sec:numerical_results}.}
    \label{tab:training_times}
    \renewcommand{\arraystretch}{1.05}
    \setlength{\tabcolsep}{4.5pt}
    \begin{tabular}{lllcc}
        \toprule
        Problem & Model & Setting & $d_z$
        & Training time (s) \\
        \midrule
        1D1V    & VPFP & $\varepsilon=1$, $T=1$        & -- &  243.2 \\
        1D1V    & BGK--Poisson  & $\varepsilon=1$, $T=1$        & -- & 245.0 \\
        2D2V    & VPFP & $\varepsilon=10^{-3}$, $N=32$ & -- & 439.4 \\
        3D3V    & BGK--Poisson  & $\varepsilon=1$               & -- &  673.1 \\
        3D3V    & BGK--Poisson  & $\varepsilon=10^{-3}$         & -- &  730.2 \\
        3D3V UQ & VPFP & $\varepsilon=10^{-4}$         & 3  & 1352.3 \\
        3D3V UQ & VPFP & $\varepsilon=10^{-4}$         & 5  & 1524.8 \\
        \bottomrule
    \end{tabular}
\end{table}

\FloatBarrier
\section{Conclusion and discussion}

In this paper, an asymptotic-preserving tensor neural network (APTNN) method is proposed for high-dimensional kinetic equations with uncertainty in the high-field regime, particularly for the Vlasov--Poisson--Fokker--Planck system and the BGK-Poisson system. To alleviate the curse of dimensionality,  we adopt TNN to approximate the solutions as sums of products of one-dimensional functions and apply dimension-wise integration to evaluate the loss, avoiding costly computations and storage on full tensor grids. To improve prediction accuracy in the high-field limit, we further incorporate the first velocity-moment equation into the baseline AP loss to explicitly constrain the flux closure. The enhanced AP loss can still be evaluated dimension-wise. Numerical results demonstrate that the APTNN method significantly reduces training time and memory usage while improving prediction accuracy in the high-field regime. Therefore, the proposed method integrates accurate prediction in the macroscopic limit with efficient high-dimensional computation, achieving large-scale kinetic simulations through economical GPU usage.

The proposed method merits investigation in several directions. First, further reducing memory consumption is important for extending APTNN to a wider range of physical problems. Kinetic solutions do not always admit low-rank approximations. For example, nonlinear plasma instabilities \cite{einkemmer2025review} and the transport of discontinuities \cite{guo2022low} may require larger ranks to maintain accuracy. The coupling between the distribution function and the electric field introduces multiplicative rank dependence into the number of separated residual terms. As shown in Corollary~\ref{cor:aptnn_rank_complexity}, Gram-matrix storage for loss evaluation grows quadratically with this number, making memory a potential bottleneck at larger ranks. Promising ways to reduce memory usage include blockwise Gram-matrix evaluation \cite{CharlierEtAl2021} and recomputation during backpropagation \cite{ChenEtAl2016}.  Second, the quantitative relationship between the separation ranks and approximation errors of TNNs requires further investigation. The current approximation results  guarantee arbitrary accuracy but do not quantify the ranks needed to achieve it. Developing rank-dependent error estimates would guide adaptive rank selection and help balance accuracy and computational cost.

\appendix
\titleformat{\section}
  {\fontsize{13pt}{15pt}\selectfont\bfseries}{\appendixname~\thesection}{1em}{}
\section{Proof of Theorem~\ref{thm:whole_velocity_weighted_tnn}}\label{app:whole_velocity_tnn_proof}

\begin{proof}
Let $\delta>0$. We first approximate the distribution $f$ by a function with compact velocity support. Specifically, choose a cutoff function $\chi\in C_c^\infty(\mathbb R)$ with $0\leq\chi\leq1$, $\chi=1$ on $[-1,1]$, and $\operatorname{supp}\chi\subset[-2,2]$. For $R\geq1$, define $\chi_R(v)=\prod_{j=1}^3\chi(v_j/R)$ and $f_R=\chi_R f$. Then $f_R=f$ on $v\in[-R,R]^3$ and $f_R$ vanishes outside $[-2R,2R]^3$ in velocity.

To show that $f_R$ converges to $f$ in $H_s^m(\mathcal D)$, we estimate the mixed derivatives of the truncation error $f-f_R=(1-\chi_R)f$. For any $\alpha,\beta,\gamma$ satisfying $\alpha+|\beta|+|\gamma|\leq m$, the Leibniz formula gives
\begin{equation}\label{eq:whole_uat_cutoff_derivatives}
\partial_t^\alpha\partial_x^\beta\partial_v^\gamma(f-f_R)=(1-\chi_R)\partial_t^\alpha\partial_x^\beta\partial_v^\gamma f-\sum_{\substack{\nu\leq\gamma\\\nu\neq0}}
\binom{\gamma}{\nu}(\partial_v^\nu\chi_R)
\partial_t^\alpha\partial_x^\beta\partial_v^{\gamma-\nu}f,
\end{equation}
where $\nu\in\mathbb N_0^3$ and $\nu\leq\gamma$ means $\nu_j\leq\gamma_j$ for $j=1,2,3$. In the second term, we use the chain rule to obtain
\[
\partial_v^\nu\chi_R(v)=R^{-|\nu|}
\prod_{j=1}^3\chi^{(\nu_j)}(v_j/R),
\text{ and }
\|\partial_v^\nu\chi_R\|_{L^\infty(\mathbb R^3)}\leq C_\nu R^{-|\nu|}\leq C_\nu,
\]
where $C_\nu=\prod_{j=1}^3 \|\chi^{(\nu_j)}\|_{L^\infty(\mathbb R)}$. Then $\alpha+|\beta|+|\gamma-\nu|=\alpha+|\beta|+|\gamma|-|\nu|\leq m$ implies that only derivatives of $f$ up to order $m$ are involved in the second term. Moreover, since $\chi_R=1$ on $(-R,R)^3$, both the coefficients $1-\chi_R$ and $\partial_v^\nu\chi_R$ for $\nu\neq0$ vanish when $|v|<R$. Therefore, using \eqref{eq:whole_uat_cutoff_derivatives} and the Cauchy--Schwarz inequality, we obtain the estimate
\begin{equation}\label{eq:whole_uat_tail_bound}
\begin{aligned}
\|f-f_R\|_{H_s^m(\mathcal D)}^2
&\leq C\sum_{\alpha+|\beta|+|\gamma|\leq m}\int_{\mathcal D\cap\{|v|\geq R\}}\langle v\rangle^{2s}\left|\partial_t^\alpha\partial_x^\beta\partial_v^\gamma f\right|^2\,dt\,dx\,dv,
\end{aligned}
\end{equation}
where $C$ is independent of $R$. Since $f\in H_s^m(\mathcal D)$, the right-hand side of \eqref{eq:whole_uat_tail_bound} tends to zero as $R\to\infty$. Thus, we can choose $R$ sufficiently large such that
\begin{equation}\label{eq:whole_uat_truncation_error}
\|f-f_R\|_{H_s^m(\mathcal D)}<\frac{\delta}{2}.
\end{equation}

We next construct a TNN approximation of $f_R$ on the whole velocity space, starting with an approximation on a bounded domain. Choose $V>2R$ and define the bounded domain $K_V=(0,T)\times\Omega\times(-V,V)^3$. Since $\langle v\rangle^s\geq1$, one has  $f_R\in H^m(K_V)$ with $\|f_R\|_{H^m(K_V)}\leq\|f_R\|_{H_s^m(\mathcal D)}<\infty$. By Theorem~\ref{thm:bounded_product_tnn}, there exists a TNN $U_\theta$ that approximates $f_R$ arbitrarily closely in $H^m(K_V)$ and has the form
\begin{equation}\label{eq:whole_uat_bounded_network}
U_\theta(t,x,v)=\sum_{r=1}^p A_{r,\theta}(t)
\prod_{i=1}^d B_{r,i,\theta}(x_i)
\prod_{j=1}^3 D_{r,j,\theta}(v_j),
\end{equation}
where $A_{r,\theta}$, $B_{r,i,\theta}$ and $D_{r,j,\theta}$ are the $r$th output components of the subnetworks for $t$, $x_i$ and $v_j$, respectively. Then, to obtain an approximation on the whole velocity space, we multiply the velocity factors $D_{r,j,\theta}$ by a smooth cutoff function and extend them by zero. Specifically, choose $\eta\in C_c^\infty(-V,V)$ with $\eta=1$ on $[-2R,2R]$ and define $C_{r,j,\theta}(v_j)=\eta(v_j)D_{r,j,\theta}(v_j)$. Each $C_{r,j,\theta}$ is then extended by zero outside $(-V,V)$. Replacing $D_{r,j,\theta}$ with $C_{r,j,\theta}$ in $U_\theta$, we obtain
\begin{equation}\label{eq:whole_uat_separated_function}
S_\theta(t,x,v)=\sum_{r=1}^p A_{r,\theta}(t)
\prod_{i=1}^d B_{r,i,\theta}(x_i)
\prod_{j=1}^3 C_{r,j,\theta}(v_j).
\end{equation}
Both $f_R$ and $S_\theta$ vanish on $\mathcal D\setminus K_V$, so it suffices to estimate their difference on $K_V$. Since $\prod_{j=1}^3\eta(v_j)=1$ on $\operatorname{supp}_v f_R$, we have
\begin{equation}\label{eq:whole_uat_separated_error_identity}
f_R-S_\theta=\left(\prod_{j=1}^3\eta(v_j)\right)(f_R-U_\theta)\quad\text{on }K_V.
\end{equation}
Taking the $H^m(K_V)$-norm in
\eqref{eq:whole_uat_separated_error_identity}
and applying the Leibniz formula, we obtain
\begin{equation}\label{eq:whole_uat_unweighted_error}
\|f_R-S_\theta\|_{H^m(K_V)}
\leq C_{m,\eta}
\|f_R-U_\theta\|_{H^m(K_V)},
\end{equation}
where $C_{m,\eta}$ depends only on $m$ and $\eta$. Together with the bound $\langle v\rangle^s\leq(1+3V^2)^{s/2}$ on $(-V,V)^3$, \eqref{eq:whole_uat_unweighted_error} yields
\begin{equation}\label{eq:whole_uat_separated_error}
\begin{aligned}
\|f_R-S_\theta\|_{H_s^m(\mathcal D)}
&\leq (1+3V^2)^{s/2}
      \|f_R-S_\theta\|_{H^m(K_V)}\\
&\leq C_{V,\eta}\|f_R-U_\theta\|_{H^m(K_V)}
<\frac{\delta}{2},
\end{aligned}
\end{equation}
where $C_{V,\eta}=(1+3V^2)^{s/2}C_{m,\eta}$. The last inequality follows by choosing $U_\theta$ such that
$\|f_R-U_\theta\|_{H^m(K_V)}<\delta/(2C_{V,\eta})$.

Finally, combining \eqref{eq:whole_uat_truncation_error} and \eqref{eq:whole_uat_separated_error}, the triangle inequality gives
\[
\begin{aligned}
\|f-S_\theta\|_{H_s^m(\mathcal D)}
&\leq \|f-f_R\|_{H_s^m(\mathcal D)}
+\|f_R-S_\theta\|_{H_s^m(\mathcal D)}\\
&<\frac{\delta}{2}+\frac{\delta}{2}=\delta.
\end{aligned}
\]

\end{proof}

\section*{Acknowledgments}
SJ is partially supported by the NSFC grant No. 12031013, and the Shanghai Commission of S\&T grant No. 23JC1402300.
The computations in this paper were run on the Siyuan-1 cluster supported by the Center for High Performance Computing at Shanghai Jiao Tong University.

\bibliographystyle{siam}
\bibliography{ref}
\end{document}